\documentclass[letterpaper,11pt]{article}
\usepackage{fullpage,fancyhdr,graphicx}
\usepackage{amsmath,amssymb,amsthm,stmaryrd}
\usepackage{tikz,pgfplots,cite}
\usepackage[colorlinks, linkcolor=blue, citecolor=violet,
urlcolor=magenta, hypertexnames=false, pageanchor=true, unicode]{hyperref}

\usetikzlibrary{shapes.multipart,positioning}
\pgfplotsset{compat=1.16}

\allowdisplaybreaks[4]

\newcommand{\BEAS}{\begin{eqnarray*}}
\newcommand{\EEAS}{\end{eqnarray*}}
\newcommand{\BEA}{\begin{eqnarray}}
\newcommand{\EEA}{\end{eqnarray}}
\newcommand{\BEQ}{\begin{equation}}
\newcommand{\EEQ}{\end{equation}}
\newcommand{\BIT}{\begin{itemize}}
\newcommand{\EIT}{\end{itemize}}
\newcommand{\ie}{{\it i.e.}}

\newcommand{\argmin}{\mathop{\rm argmin}}

\newcommand{\prox}{\mathrm{prox}}
\newcommand{\dom}{\mathop{\bf dom}}
\newcommand{\intr}{\mathop{\bf int}}
\newcommand{\inprod}[2]{\langle#1, #2\rangle}

\newcommand{\ones}{\mathbf 1}
\newcommand{\reals}{{\mbox{\bf R}}}

\newcommand{\avg}{\mathrm{avg}}
\newcommand{\primal}{\mathrm p}
\newcommand{\dual}{\mathrm d}
\newcommand{\cL}{\mathcal L}
\newcommand{\red}[1]{{\color{red}{#1}}}

\makeatletter
\renewcommand*\env@matrix[1][c]{\hskip -\arraycolsep
  \let\@ifnextchar\new@ifnextchar
  \array{*\c@MaxMatrixCols #1}}
\makeatother

\title{Bregman three-operator splitting methods}
\author{Xin Jiang%
\thanks{Department of Electrical and Computer Engineering, UCLA.
Email: \textsf{jiangxjames@ucla.edu}, \textsf{vandenbe@ucla.edu}.}
\and Lieven Vandenberghe\footnotemark[1]}
\date{October 3, 2022}

\begin{document}
\maketitle

\begin{abstract}
The paper presents primal--dual proximal splitting methods for
convex optimization, in which generalized Bregman distances are used
to define the primal and dual proximal update steps.
The methods extend the primal and dual Condat--V\~u algorithms
and the primal--dual three-operator (PD3O) algorithm.
The Bregman extensions of the Condat--V\~u algorithms 
are derived from the Bregman proximal point method applied to
a monotone inclusion problem.
Based on this interpretation, a unified framework for the convergence 
analysis of the two methods is presented.
We also introduce a line search procedure for stepsize selection
in the Bregman dual Condat--V\~u algorithm applied to 
equality-constrained problems.
Finally, we propose a Bregman extension of PD3O and analyze its 
convergence.
\end{abstract}

\section{Introduction}

We discuss proximal splitting methods for optimization problems 
in the form
\begin{equation} \label{e-prob}
  \mbox{minimize} \quad f(x) + g(Ax) + h(x),
\end{equation}
where $f$, $g$, and $h$ are convex functions, and $h$ is differentiable.
This general problem covers a wide variety of applications in 
machine learning, signal and image processing, operations research,
control, and other fields~\cite{CP11b,PaB:13,KP15,CP16b}.
In this paper, we consider proximal splitting methods based on Bregman 
distances for solving~\eqref{e-prob}
and some interesting special cases of~\eqref{e-prob}.

Recently, several primal--dual first-order methods have been proposed for 
the three-term problem~(\ref{e-prob}):
the Condat--V\~u algorithm~\cite{Con:13,Vu:13,YETA22},
the primal--dual three-operator (PD3O) algorithm \cite{Yan18}, and
the primal--dual Davis--Yin (PDDY) algorithm~\cite{SCMR20}.
Algorithms for some special cases of~\eqref{e-prob} are also of interest.
These include the Chambolle--Pock algorithm, also known as the
primal--dual hybrid gradient (PDHG) method~\cite{ChP:11,ChP:16} 
(when $h=0$),
the Loris--Verhoeven algorithm~\cite{LV11,CHZ13,DST15} (when $f=0$), 
the proximal gradient algorithm (when $g=0$), and
the Davis--Yin splitting algorithm~\cite{DY15} (when $A=I$).
All these methods handle the nonsmooth functions $f$ and $g$
via the standard Euclidean proximal operator.
 
To further improve the efficiency of proximal algorithms, 
proximal operators based on generalized Bregman distances have 
been proposed and incorporated in many methods
\cite{ChT:93,Eck:93,Gul:94,AuT:06,Tse:08,BBT17,BSTV18,LFN18,Teboulle18}. 
Bregman distances offer two potential benefits.
First, the Bregman distance can help build a more accurate local
optimization model around the current iterate.
This is often interpreted as a form of preconditioning.  For example,
diagonal or quadratic preconditioning~\cite{PC11,JLLO19,LXY21}
has been shown to improve the practical convergence of PDHG,
as well as the accuracy of the computed solution~\cite{ADH+:21}.
As a second benefit, a Bregman proximal operator of a 
function may be easier to compute than 
the standard Euclidean proximal operator, and therefore
reduce the complexity per iteration of an optimization algorithm.
Recent applications of this kind include optimal transport 
problems~\cite{CLMW19}, optimization over nonnegative trigonometric 
polynomials \cite{ChV:18},
and sparse semidefinite programming~\cite{JV22}.

Extending standard proximal methods and their convergence analysis
to Bregman distances is not straightforward because some fundamental
properties of the Euclidean proximal operator no longer hold
for Bregman proximal operators.
An example is the Moreau decomposition
which relates the (Euclidean) proximal operators of a closed convex 
function and its conjugate~\cite{Mor:65}.
Another example is the simple relation between
the proximal operators of a function $g$ and 
the composition with a linear function $g(Ax)$ 
when $AA^T$ is a multiple of the identity;
see, \emph{e.g.},~\cite{CP11b,Beck17}.
This composition rule is used in~\cite{OV20}
to establish the equivalence between some well-known first-order proximal 
methods for problem~\eqref{e-prob} with $A=I$ and with general $A$.

The purpose of this paper is to present new Bregman extensions and
convergence results for the Condat--V\~u and PD3O algorithms.
The main contributions are as follows.
\begin{itemize}
\item The Condat--V\~u algorithm~\cite{Con:13,Vu:13} exists in a primal
and a dual variant.
We discuss extensions of the two algorithms
that use Bregman proximal operators in the primal and dual updates.
The Bregman primal Condat--V\~u algorithm first appeared
in~\cite[Algorithm 1]{ChP:16}, and is also a special case 
of the algorithm proposed in~\cite{YHA21}
for a more general convex--concave saddle point problem.
We give a new derivation of this method and its dual variant,
by applying the Bregman proximal point method
to the primal--dual optimality conditions.
Based on the interpretation, we provide a unified framework
for the convergence analysis of the two variants,
and show an $O(1/k)$ ergodic convergence rate,
which is consistent with previous results for Euclidean proximal operators
in~\cite{Con:13,Vu:13} and Bregman proximal operators in~\cite{ChP:16}.
We also give a convergence result for the primal and dual iterates. 

\item We propose an easily implemented backtracking line search technique 
for selecting stepsizes in the Bregman dual Condat--V\~u
algorithm for problems with equality constraints.
The proposed backtracking procedure is similar to the technique
in~\cite{MaP:18} for the special setting of PDHG with Euclidean proximal 
operators, but has some important differences even in this special case.
We give a detailed analysis of the algorithm
with line search and recover the $O(1/k)$ ergodic rate of convergence
for related algorithms in~\cite{MaP:18,JV22}.

\item We propose a Bregman extension for PD3O
and establish an ergodic convergence result.
\end{itemize}
The paper is organized as follows.
Section~\ref{s-duality} gives a precise statement of the 
problem~(\ref{e-prob}), and reviews the duality theory
that will be used in the rest of the paper. 
In Section~\ref{s-overview} we review some well-known first-order
proximal methods and establish connections between them.
Section~\ref{s-bregman} provides some necessary background on 
Bregman distances.
In Section~\ref{s-bcv} we discuss the Bregman primal and dual 
Condat--V\~u algorithms
and analyze their convergence.
The line search technique and its convergence are discussed
in Section~\ref{s-ls}.
In Section~\ref{s-bpd3o} we extend PD3O to a Bregman proximal method
and analyze its convergence.
Section~\ref{s-exp} contains results of a numerical experiment.

\section{Duality theory and merit functions} \label{s-duality}

This section summarizes the facts from convex duality theory  
that underlie the primal--dual methods discussed in the paper.
We also describe primal--dual merit functions that will be used in
the convergence analysis.

We use the notation $\inprod{x}{y} = x^Ty$ for the
standard inner product of vectors $x$ and $y$, and
$\|x\| = \inprod{x}{x}^{1/2}$ for the Euclidean norm of a vector
$x$.  Other norms will be distinguished by a subscript.

\subsection{Problem formulation}

In~\eqref{e-prob} the vector $x$ is an $n$-vector and
$A$ is an $m\times n$ matrix.  The functions $f$, $g$, $h$
are closed and convex, and $h$ is differentiable, \ie,  
\[
h(x) \geq h(x') + \inprod{\nabla h(x')}{x-x'} \quad
\mbox{for all $x,x' \in \dom h$},
\]
where $\dom h$ is an open convex set.  We assume that $f+h$ and $g$
are proper, \ie, have nonempty domains.

An important example of~\eqref{e-prob} is $g = \delta_C$,
the indicator function of a closed convex set $C$.
With $g=\delta_C$, the problem is equivalent to 
\[
 \begin{array}{ll}
 \mbox{minimize} & f(x) + h(x) \\
 \mbox{subject to} & Ax \in C.
 \end{array}
\]
For $C =\{b\}$ the constraints are a set of linear equations $Ax=b$. 
This special case actually covers all applications of the more general
problem~\eqref{e-prob}, since~\eqref{e-prob} can be reformulated as
\[
 \begin{array}{ll}
 \mbox{minimize} & f(x) + g(y) + h(x) \\
 \mbox{subject to} & Ax=y,
 \end{array}
\]
at the expense of increasing the problem size
by introducing a splitting variable $y$.

\subsection{Dual problem and optimality conditions}

The dual of problem~\eqref{e-prob} is
\BEQ \label{e-prob-dual}
 \mbox{maximize}\quad  {-(f+h)^*(-A^Tz) - g^*(z)},
\EEQ
where $(f+h)^*$ and $g^*$ are the conjugates of $f+h$ and $g$:
\[
 (f+h)^*(w) = \sup_x{(\inprod{w}{x} - f(x) - h(x))}, \qquad
 g^*(z) = \sup_y{(\inprod{z}{y} - g(y))}.
\]
The conjugate $(f+h)^*$ is the infimal convolution 
$f^*$ and $h^*$, denoted by $f^* \boxempty h^*$:
\[
(f^* \boxempty h^*) (z) = \inf_w{((f^*(w) + h^*(z-w))}.
\]

The primal--dual optimality conditions for~\eqref{e-prob}
and~\eqref{e-prob-dual} are
\[
 0 \in \partial f(x) + \nabla h(x) + A^T z, \qquad
 0 \in \partial g^*(z) - Ax.
\]
Here $\partial f$ and $\partial g^*$ are the subdifferentials of
$f$ and  $g^*$.  We often write the optimality conditions as
\begin{equation} \label{e-opt-cond}
  0 \in \begin{bmatrix} 0 & A^T \\ -A & 0 \end{bmatrix}
  \begin{bmatrix} x \\ z \end{bmatrix} + \begin{bmatrix}
    \partial f(x)+\nabla h(x) \\ \partial g^\ast(z)
\end{bmatrix}.
\end{equation}
Throughout the paper, we assume that the optimality 
conditions~\eqref{e-opt-cond} are solvable.

We will refer to the convex--concave function
\[
\cL(x,z) = f(x)+h(x)+ \inprod{z}{Ax} -g^\ast(z)
\]
as the \emph{Lagrangian} of~\eqref{e-prob}.
We follow the convention that $\cL(x,z) = +\infty$ if $x\not\in \dom(f+h)$
and $\cL(x,z) = -\infty$ if $x\in\dom(f+h)$ and $z\not\in \dom g^*$.
The objective functions in~\eqref{e-prob} and
the dual problem~\eqref{e-prob-dual} can be expressed as
\[
\sup_z \cL(x,z) = f(x)+h(x)+g(Ax), \qquad
\inf_x \cL(x,z)  = -(f+h)^\ast(-A^Tz)-g^\ast(z).
\]
Solutions $x^\star$, $z^\star$ of the optimality 
conditions~\eqref{e-opt-cond} form a saddle-point of $\cL$, \ie, satisfy
\BEQ \label{e-gap}
\inf_x \sup_z \cL(x, z) = \sup_z \cL(x^\star, z) = 
\cL(x^\star, z^\star) = 
\inf_x \cL(x,z^\star) = \sup_z \inf_x \cL(x,z).
\EEQ
In particular, $\cL(x^\star, z^\star)$ is the optimal value 
of~\eqref{e-prob} and~\eqref{e-prob-dual}.

\subsection{Merit functions}

The algorithms discussed in this paper generate primal and dual iterates 
and approximate solutions $x$, $z$ with 
$x \in \dom (f+h)$ and $z\in \dom g^*$.
The feasibility conditions $Ax \in \dom g$
and $-A^Tz \in \dom (f+h)^*$ are not necessarily satisfied.
Hence the duality gap
\BEQ \label{e-gap2}
\sup_{z'} \cL(x, z') - \inf_{x'} \cL(x', z) 
 =  f(x) + h(x) + g(Ax) + (f+h)^*(-A^Tz) + g^*(z)
\EEQ
may not always be useful as a merit function to measure convergence. 

If we add constraints $x'\in X$ and $z'\in Z$ to the optimization
problems on the left-hand side of~\eqref{e-gap2},
where $X$ and $Z$ are compact convex sets, we obtain a function
\BEQ \label{e-merit}
\eta(x,z) = \sup_{z'\in Z} \cL(x, z') - \inf_{x'\in X} \cL(x', z) 
\EEQ
defined for all $x\in\dom(f+h)$ and $z\in\dom g^*$.
This follows from the fact that the functions $f+h+\delta_X$ and
$g^*+\delta_Z$ are closed and co-finite, so their conjugates
have full domain \cite[Corollary 13.3.1]{Roc:70}. 
If $\eta(x,z)$ is easily computed, and $\eta(x,z) \geq 0$ for all 
$x\in \dom(f+h)$ and $z\in\dom g^*$ with equality only if 
$x$ and $z$ are optimal, then the function $\eta$
can serve as a merit function in primal--dual algorithms
for problem~\eqref{e-prob}.

If $\dom(f+h)$ and $\dom g^*$ are bounded,
then $X$ and $Z$ can be chosen to contain $\dom (f+h)$ and $\dom g^*$.
Then the constraints in~\eqref{e-merit} are redundant and $\eta(x,z)$
is the duality gap~\eqref{e-gap2}.
Boundedness of $\dom(f+h)$ and $\dom g^*$ is a common
assumption in the literature on primal--dual first-order methods.

A weaker assumption is that~\eqref{e-prob}
has an optimal solution $x^\star \in \intr(X)$
and~\eqref{e-prob-dual} has an optimal solution $z^\star\in \intr(Z)$.
Then $\eta(x,z) \geq 0$
for all $x \in \dom(f+h)$ and $z \in \dom g^*$, with equality
$\eta(x,z) = 0$ only if $x,z$ are optimal for~\eqref{e-prob}
and~\eqref{e-prob-dual}.
To see this, we first express the two terms in~\eqref{e-merit}~as
\begin{align*}
\sup_{z'\in Z} \cL(x,z') &= f(x) + h(x) + (g \boxempty \sigma_Z) (Ax), \\
\inf_{x'\in X} \cL(x',z) &= -g^*(z) - ( (f+h)^* \boxempty \sigma_X)(-A^Tz),
\end{align*}
where $\sigma_X = \delta_X^*$ and $\sigma_Z(v) = \delta_Z^*$
are the support functions of $X$ and $Z$, respectively.
Consider the problem of minimizing $\eta(x,z)$.
By expanding the infimal convolutions in the expressions for the
two terms of $\eta$, this convex optimization problem can be formulated as
\BEQ \label{e-min-gap-prob}
\begin{array}{l@{\hskip 1em}l}
\mbox{minimize} &
f(x) + h(x) + g(y) + \sigma_Z(Ax-y) \\
& \mbox{} + g^*(z) + (f+h)^*(w) + \sigma_X(-A^Tz-w),
\end{array}
\EEQ
with variables $x, y, z, w$.
The dual of this problem is
\BEQ \label{e-min-gap-prob-dual}
\begin{array}{l@{\hskip 1em}l}
\mbox{maximize} & -f(\bar x) - h(\bar x) - g(A\bar x)
 -g^*(\bar z) - (f+h)^*(-A^T\bar z) \\
\mbox{subject to} & \bar x \in X, \; \bar z  \in Z,
\end{array}
\EEQ
with variables $\bar x, \bar z$.  The optimality conditions
for~(\ref{e-min-gap-prob}) and~(\ref{e-min-gap-prob-dual}) include the
conditions
$Ax-y \in N_Z(\bar z)$ and $-A^Tz - w \in N_X(\bar x)$,
where $N_X(\bar x) = \partial \delta_X(\bar x)$ is the normal cone
to $X$ at $\bar x$, and
$N_Z(\bar z) = \partial \delta_Z(\bar z)$ the normal cone
to $Z$ at $\bar z$.
By assumption, there exist points $x^\star \in \intr(X)$ and
$z^\star \in \intr(Z)$ that are optimal for the original
problem~(\ref{e-prob}) and  its dual~(\ref{e-prob-dual}).
It can be verified that
$(x,y,z,w) = (x^\star, Ax^\star, z^\star, -A^Tz^\star)$,
$(\bar x, \bar z) = (x^\star, z^\star)$
are optimal for~(\ref{e-min-gap-prob}) and~(\ref{e-min-gap-prob-dual}),
and that $\eta(x^\star,z^\star) = 0$.
Now let $(\hat x,\hat z)$ be any other minimizer of $\eta$, \ie,
$\eta(\hat x, \hat z) = 0$.
Then $\hat x, \hat z$ and the corresponding minimizers $\hat y, \hat w$
in~(\ref{e-min-gap-prob}), must satisfy the optimality conditions
with the optimal dual variables
$\bar x= x^\star$, $\bar z = z^\star$. In particular,
$A\hat x -\hat y \in N_Z(z^\star) = \{0\}$
and $-A^T\hat z -\hat w \in N_X(x^\star) = \{0\}$.
The objective value of~(\ref{e-min-gap-prob}) at this point then reduces
to $0 = f(\hat x) + h(\hat x) + g(A\hat x)
 + g^*(\hat z) + (f+h)^*(-A^T\hat w)$,  the duality gap associated
with the original problem and its dual.
This shows that $\eta(\hat x, \hat z) =0$ implies that $\hat x, \hat z$
are optimal for problem~(\ref{e-prob}) and~(\ref{e-prob-dual}).

Consider for example the primal and dual pair
\[
\begin{array}[t]{ll}
  \mbox{minimize} & f(x) + h(x) \\
  \mbox{subject to} & Ax = b
\end{array} \qquad \quad
\begin{array}[t]{ll}
  \mbox{maximize} & -b^Tz - (f+h)^*(-A^Tz).
\end{array} 
\]
Here $g= \delta_{\{b\}}$.
If we take $Z = \{z \mid \|z\| \leq \gamma\}$,
then $\sigma_Z(y)=\gamma \|y\|$ and 
$(g \boxempty \sigma_Z)(Ax) = \gamma \|Ax-b\|$.
If in addition $\dom (f+h)$ is bounded and
we take $X \supseteq \dom(f+h)$, then
\[
\eta(x,z) = f(x) + h(x) + \gamma\|Ax-b\|
+ b^Tz + (f+h)^*(-A^Tz)
\]
with domain $\dom(f+h) \times \reals^m$.
The first three terms are the primal objective augmented with an
exact penalty for the constraint $Ax=b$.

As another example, consider 
\[
\begin{array}[t]{ll}
  \mbox{minimize} & \|x\|_1 \\
  \mbox{subject to} & Ax\leq b
\end{array} \qquad\quad
\begin{array}[t]{ll}
  \mbox{maximize} & -b^T z\\
  \mbox{subject to} & \|A^Tz\|_\infty \leq 1 \\
  & z \geq 0.
\end{array}
\]
This is an example of~\eqref{e-prob} with $f(x)=\|x\|_1$, $h(x)=0$,
and $g$ the indicator function of $\{y \mid y \leq b\}$.
The domains $\dom f$ and $\dom g^*$ are unbounded.
If we choose $X = \{ x \mid \|x\|_\infty \leq \kappa\}$
and $Z = \{z \mid \|z\|_\infty \leq \lambda\}$, then
\[
\sigma_X(w) = \kappa \|w\|_1, \qquad
(f^* \boxempty \sigma_X)(w)
 = \kappa \sum_{i=1}^n \max \{0, |w_i|-1\}
\]
and
\[
\sigma_Z(y) =  \lambda \sum_{i=1}^m \max\{0, y_i\}, \qquad
(g \boxempty \sigma_Z)(y) =  \lambda \sum_{i=1}^m \max\{0, y_i-b_i\}.
\]
Hence, for this example the merit function~\eqref{e-merit} is
\[
\eta(x,z) 
= \|x\|_1 + \lambda \sum_{i=1}^m \max\{0, (Ax-b)_i\}
 + b^T z  + \kappa \sum_{i=1}^n  \max\{0, |(A^Tz)_i|-1\}
\]
with domain $\reals^n \times \reals^m_+$.
The second term is an exact penalty for the primal constraint $Ax\leq b$.
The last term is an exact penalty for the dual constraint 
$\|A^Tz\|_\infty\leq 1$.

\section{First-order proximal algorithms: survey and connections}
\label{s-overview}

In this section, we discuss several first-order proximal algorithms and 
their connections.
We start with four three-operator splitting algorithms
for problem~\eqref{e-prob}:
the primal and dual variants of the
Condat--V\~u algorithm~\cite{Con:13,Vu:13},
the primal--dual three-operator (PD3O) algorithm~\cite{Yan18},
and the primal--dual Davis--Yin (PDDY) algorithm~\cite{SCMR20}.
For each of the four algorithms, we make connections 
with other first-order proximal algorithms,
using reduction (\ie, setting some parts in~\eqref{e-prob} to zero)
and the ``completion'' reformulation (based on extending $A$ to
a matrix with orthogonal rows and equal row norms)~\cite{OV20}.
We focus on the formal connections between algorithms.
The connections do not necessarily provide the best approach
for convergence analysis or the best known convergence results.

The \emph{proximal operator} or \emph{proximal mapping} of a closed convex 
function $f \colon \reals^n \rightarrow \reals$ is defined~as
\BEQ \label{e-prox}
\prox_f(y) = \argmin_x{\big(f(x)+\frac{1}{2} \|x-y\|^2\big)}.
\EEQ
If $f$ is closed and convex, the minimizer in the definition exists 
and is unique for all $y$ \cite{Mor:65}.
We will call~\eqref{e-prox} the \emph{standard} or the \emph{Euclidean
proximal operator} when we need to distinguish it from Bregman
proximal operators defined in Section~\ref{s-bregman}.

\subsection{Condat--V\~u three-operator splitting algorithm}

\begin{figure}
\centering
\begin{tikzpicture}[>=latex,font=\small,
                    every text node part/.style={align=center}]
\node (DY) [rectangle,draw] 
  {reduced primal \\ Condat--V\~u};
\node (DR) [rectangle,draw,right=2.0cm of DY] 
  {(primal) \\ Douglas--Rachford};
\node (PG) [rectangle,draw, left=2.0cm of DY] 
  {reduced Loris--Verhoeven \\ with shift~\eqref{e-spg}};
\node (PD3O) [rectangle,draw,above=2.05cm of DY,fill=gray!30] 
  {\red{(primal)} \\ \red{Condat--V\~u} \eqref{e-cv}};
\node (PDHG) [rectangle,draw,above=2.20cm of DR] 
  {(primal) PDHG};
\node (LV) [rectangle,draw,above=1.98cm of PG]
  {Loris--Verhoeven \\ with shift~\eqref{e-slv}};
\node (PG0) [rectangle,draw,above=1.0cm of PD3O] {proximal gradient};
\path (DY.north) -- (DY.north east) coordinate[pos=0.36] (DYe);
\path (DY.north) -- (DY.north west) coordinate[pos=0.36] (DYw);
\path (DR.north) -- (DR.north east) coordinate[pos=0.30] (DRe);
\path (DR.north) -- (DR.north west) coordinate[pos=0.30] (DRw);
\path (PG.north) -- (PG.north east) coordinate[pos=0.5] (PGe);
\path (PG.north) -- (PG.north west) coordinate[pos=0.5] (PGw);
\path (PD3O.south) -- (PD3O.south east) coordinate[pos=0.33] (PD3Oe);
\path (PD3O.south) -- (PD3O.south west) coordinate[pos=0.35] (PD3Ow);
\path (PDHG.south) -- (PDHG.south east) coordinate[pos=0.33] (PDHGe);
\path (PDHG.south) -- (PDHG.south west) coordinate[pos=0.33] (PDHGw);
\path (LV.south) -- (LV.south east) coordinate[pos=0.36] (LVe);
\path (LV.south) -- (LV.south west) coordinate[pos=0.36] (LVw);
\draw[-latex] (DY)   -- node[midway, above] {$h=0$} (DR);
\draw[-latex] (DY)   -- node[midway, above] {$f=0$} (PG);
\draw[-latex] (PD3O) -- node[midway, above] {$h=0$} (PDHG);
\draw[-latex] (PD3O) -- node[midway, above] {$f=0$} (LV);
\draw[-latex] (DYw)  -- node[midway, above, sloped] {completion}
  (PD3Ow);
\draw[latex-] (DYe)  -- node[midway, right]  {$A=I$} (PD3Oe);
\draw[-latex] (DRw)  -- node[midway, above, sloped] {completion} 
  (PDHGw);
\draw[latex-] (DRe)  -- node[midway, right]  {$A=I$} (PDHGe);
\draw[latex-] (PG)   -- node[midway, right]  {$A=I$} (LV);
\draw[latex-] (PG0)  -- node[midway, right]  {$g=0$} (PD3O);
\end{tikzpicture}
\caption{Proximal methods derived from primal Condat--V\~u algorithm.}
\label{t-sum-cv}
\end{figure}
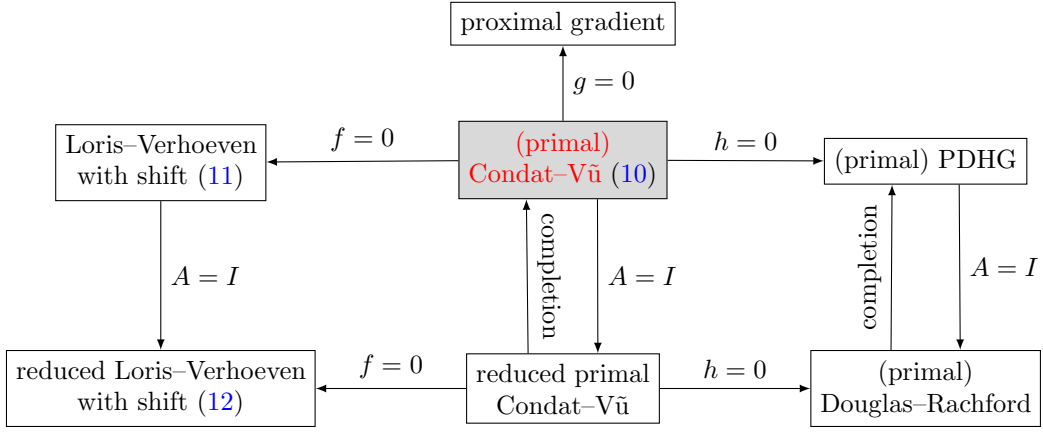

We start with the (primal) Condat--V\~u three-operator splitting 
algorithm, which was proposed independently
by Condat~\cite{Con:13} and V\~u~\cite{Vu:13},
\begin{subequations} \label{e-cv}
\begin{align}
  x^{(k+1)} &= \prox_{\tau f}
    \big(x^{(k)}-\tau (A^Tz^{(k)}+\nabla h(x^{(k)})) \big) \\
  z^{(k+1)} &= \prox_{\sigma g^\ast}
    \big(z^{(k)}+\sigma A(2x^{(k+1)}-x^{(k)}) \big).
\end{align}
\end{subequations}
The stepsizes $\sigma$ and $\tau$ must satisfy
\[
\sigma \tau \|A\|_2^2 + \tau L \leq 1,
\]
where $\|A\|_2$ is the spectral norm of $A$, and~$L$ is
the Lipschitz constant of $\nabla h$ with respect to the Euclidean norm.
Many other first-order proximal algorithms can be viewed as
special cases of~\eqref{e-cv},
and their connections are summarized in Figure~\ref{t-sum-cv}.
When $h=0$, algorithm~\eqref{e-cv} reduces to the (primal)
primal--dual hybrid gradient (PDHG) method~\cite{PCBC09,ChP:11,ChP:16},
or PDHGMu in~\cite{EZC:10}.
When $g=0$ in~\eqref{e-cv} (and assuming $z^{(0)}=0$),
we obtain the proximal gradient algorithm.
When $f=0$, we obtain a variant of the Loris--Verhoeven algorithm
\cite{LV11,CHZ13,DST15},
\begin{subequations} \label{e-slv}
\begin{align}
  x^{(k+1)} &= x^{(k)}-\tau (A^Tz^{(k)}+\nabla h(x^{(k)})) \\
  z^{(k+1)} &= \prox_{\sigma g^\ast}
    \big((I-\sigma \tau AA^T) z^{(k)}
    +\sigma A(x^{(k+1)}-\tau \nabla h(x^{(k)})) \big).
\end{align}
\end{subequations}
We refer to this as \emph{Loris--Verhoeven with shift},
for reasons that will be clarified later.
Furthermore, when $A=I$ in PDHG,
we obtain the Douglas--Rachford splitting (DRS) algorithm
\cite{LiM:79,EcB:92,CoP:07}.
Conversely, the ``completion'' technique in~\cite{OV20} shows that
PDHG coincides with DRS applied to a reformulation of the problem.
Similarly, when $A=I$ in the primal Condat--V\~u algorithm~\eqref{e-cv},
we obtain a new algorithm and refer to it as
the \emph{reduced primal Condat--V\~u algorithm}.
Conversely, the reduced primal Condat--V\~u algorithm
reverts to~\eqref{e-cv} via the ``completion'' trick.
We can also set $f=0$ in the reduced Condat--V\~u algorithm
or $A=I$ in~\eqref{e-slv},
and obtain the \emph{reduced Loris--Verhoeven algorithm with shift}:
\begin{subequations} \label{e-spg}
\begin{align}
  x^{(k+1)} &= x^{(k)}-\tau(z^{(k)}+\nabla h(x^{(k)})) \\
  z^{(k+1)} &= \prox_{\sigma g^\ast}
    \big( (1-\sigma\tau)z^{(k)}+\sigma(x^{(k+1)}
    -\tau \nabla
    h(x^{(\red{k})})\big).
\end{align}
\end{subequations}
Finally, due to the absence of~$f$ in~\eqref{e-spg},
it is not clear how to apply the ``completion'' trick to~\eqref{e-spg} 
to obtain~\eqref{e-slv}.

\begin{figure}
\centering
\begin{tikzpicture}[>=latex, font=\small,
                    every text node part/.style={align=center}]
\node (DY) [rectangle,draw] {reduced dual \\ Condat--V\~u};
\node (DR) [rectangle,draw,right=2.0cm of DY] 
  {dual \\ Douglas--Rachford};
\node (PG) [rectangle,draw, left=2.0cm of DY] 
  {reduced dual Loris--Verhoeven \\ with shift~\eqref{e-spg-dual}};
\node (PD3O) [rectangle,draw,above=2.05cm of DY,fill=gray!30] 
  {\red{dual} \\ \red{Condat--V\~u} \eqref{e-cv-dual}};
\node (PDHG) [rectangle,draw,above=2.28cm of DR] {dual PDHG};
\node (LV) [rectangle,draw,above=1.98cm of PG]
  {dual Loris--Verhoeven \\ with shift~\eqref{e-slv-dual}};
\node (PG0) [rectangle,draw,above=1.0cm of PD3O] {proximal gradient};
\path (DY.north) -- (DY.north east) coordinate[pos=0.42] (DYe);
\path (DY.north) -- (DY.north west) coordinate[pos=0.42] (DYw);
\path (DR.north) -- (DR.north east) coordinate[pos=0.23] (DRe);
\path (DR.north) -- (DR.north west) coordinate[pos=0.23] (DRw);
\path (PG.north) -- (PG.north east) coordinate[pos=0.5] (PGe);
\path (PG.north) -- (PG.north west) coordinate[pos=0.5] (PGw);
\path (PD3O.south) -- (PD3O.south east) coordinate[pos=0.34] (PD3Oe);
\path (PD3O.south) -- (PD3O.south west) coordinate[pos=0.34] (PD3Ow);
\path (PDHG.south) -- (PDHG.south east) coordinate[pos=0.32] (PDHGe);
\path (PDHG.south) -- (PDHG.south west) coordinate[pos=0.32] (PDHGw);
\path (LV.south) -- (LV.south east) coordinate[pos=0.36] (LVe);
\path (LV.south) -- (LV.south west) coordinate[pos=0.36] (LVw);
\draw[-latex] (DY)   -- node[midway, above] {$h=0$} (DR);
\draw[-latex] (DY)   -- node[midway, above] {$f=0$} (PG);
\draw[-latex] (PD3O) -- node[midway, above] {$h=0$} (PDHG);
\draw[-latex] (PD3O) -- node[midway, above] {$f=0$} (LV);
\draw[-latex] (DYw)  -- node[midway, above, sloped] {completion}
  (PD3Ow);
\draw[latex-] (DYe)  -- node[midway, right]  {$A=I$} (PD3Oe);
\draw[-latex] (DRw)  -- node[midway, above, sloped] {completion} 
  (PDHGw);
\draw[latex-] (DRe)  -- node[midway, right]  {$A=I$} (PDHGe);
\draw[latex-] (PG)   -- node[midway, right]  {$A=I$} (LV);
\draw[latex-] (PG0)  -- node[midway, right]  {$g=0$} (PD3O);
\end{tikzpicture}
\caption{Proximal methods derived from dual Condat--V\~u algorithm.}
\label{t-sum-cv-dual}
\end{figure}
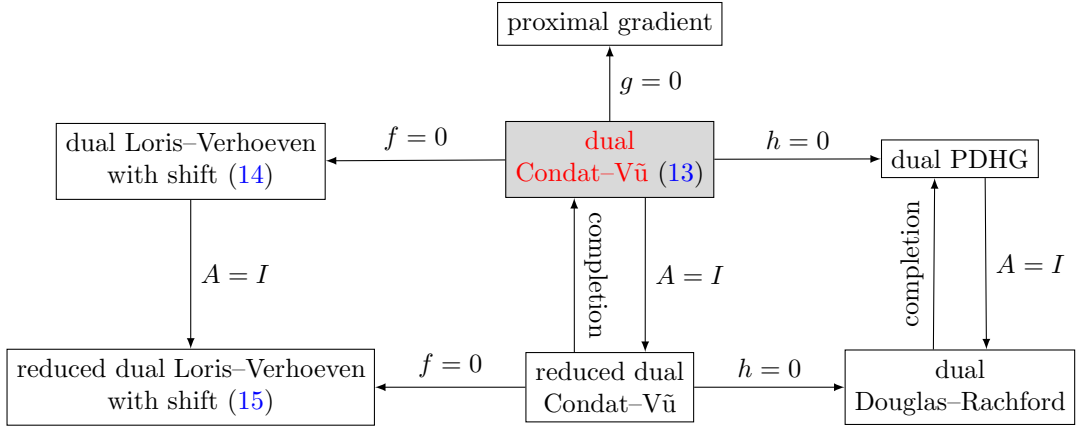

Condat~\cite{Con:13} also discusses a variant
of~\eqref{e-cv}, which we will call the 
dual Condat--V\~u algorithm:%
\begin{subequations} \label{e-cv-dual}
\begin{align}
z^{(k+1)} &= \prox_{\sigma g^\ast} (z^{(k)}+\sigma Ax^{(k)}) \\
x^{(k+1)} &= \prox_{\tau f}
  \big(x^{(k)}-\tau (A^T(2z^{(k+1)}-z^{(k)})+\nabla h(x^{(k)})) \big).
\end{align}
\end{subequations}
Figure~\ref{t-sum-cv-dual} summarizes the proximal algorithms
derived from~\eqref{e-cv-dual}.
When $h=0$, algorithm~\eqref{e-cv-dual} reduces to PDHG
applied to the dual of~\eqref{e-prob} (with $h=0$),
which is shown to be equivalent to linearized ADMM~\cite{PaB:13}
(also called Split Inexact Uzawa in~\cite{EZC:10}).
Setting $g=0$ in~\eqref{e-cv-dual} yields the proximal gradient algorithm.
When $f=0$, we obtain a new algorithm:
\begin{subequations} \label{e-slv-dual}
\begin{align}
z^{(k+1)} &= \prox_{\sigma g^\ast} (z^{(k)}+\sigma Ax^{(k)}) \\
x^{(k+1)} &= x^{(k)}-\tau (A^T(2z^{(k+1)}-z^{(k)})+\nabla h(x^{(k)})).
\end{align}
\end{subequations}
Following the previous naming convention, we call it
\emph{dual Loris--Verhoeven algorithm with shift}.
Furthermore, setting $A=I$ in~\eqref{e-cv-dual} gives
the \emph{reduced dual Condat--V\~u algorithm}.
Conversely, applying the ``completion'' trick to this reduced algorithm
recovers~\eqref{e-cv-dual}.
Similarly, setting $A=I$ in dual PDHG gives dual DRS, \ie,
DRS with $f$ and $g$ switched,
and conversely, the ``completion'' trick recovers dual PDHG
from dual DRS.
We can also set $A=I$ in~\eqref{e-slv-dual}
or $f=0$ in the reduced dual Condat--V\~u algorithm,
and obtain the \emph{reduced dual Loris--Verhoeven algorithm with shift}:
\begin{subequations} \label{e-spg-dual}
\begin{align}
z^{(k+1)} &= \prox_{\sigma g^\ast} (z^{(k)}+\sigma x^{(k)}) \\
x^{(k+1)} &= x^{(k)}-\tau (2z^{(k+1)}-z^{(k)}+\nabla h(x^{(k)})).
\end{align}
\end{subequations}

\subsection{Primal--dual three-operator (PD3O) splitting algorithm}

\begin{figure}
\centering
\begin{tikzpicture}[>=latex, font=\small,
                    every text node part/.style={align=center}]
\node (DY) [rectangle,draw] {Davis--Yin};
\node (DR) [rectangle,draw, right=2.0cm of DY] 
  {(primal) \\ Douglas--Rachford};
\node (PG) [rectangle,draw, left=2.4cm of DY] {proximal gradient};
\node (PD3O) [rectangle,draw,above=2.45cm of DY, fill=gray!30] 
  {\red{PD3O}~\eqref{e-pd3o}};
\node (PDHG) [rectangle,draw,above=2.2cm of DR] {(primal) PDHG};
\node (LV) [rectangle,draw,above=2.39cm of PG]
  {Loris--Verhoeven~\eqref{e-lv}};
\path (DY.north) -- (DY.north east) coordinate[pos=0.42] (DYe);
\path (DY.north) -- (DY.north west) coordinate[pos=0.42] (DYw);
\path (DR.north) -- (DR.north east) coordinate[pos=0.25] (DRe);
\path (DR.north) -- (DR.north west) coordinate[pos=0.25] (DRw);
\path (PG.north) -- (PG.north east) coordinate[pos=0.5] (PGe);
\path (PG.north) -- (PG.north west) coordinate[pos=0.5] (PGw);
\path (PD3O.south) -- (PD3O.south east) coordinate[pos=0.4] (PD3Oe);
\path (PD3O.south) -- (PD3O.south west) coordinate[pos=0.4] (PD3Ow);
\path (PDHG.south) -- (PDHG.south east) coordinate[pos=0.28] (PDHGe);
\path (PDHG.south) -- (PDHG.south west) coordinate[pos=0.28] (PDHGw);
\path (LV.south) -- (LV.south east) coordinate[pos=0.36] (LVe);
\path (LV.south) -- (LV.south west) coordinate[pos=0.36] (LVw);
\draw[-latex] (DY)   -- node[midway, above] {$h=0$} (DR);
\draw[-latex] (DY)   -- node[midway, above] {$f=0$} (PG);
\draw[-latex] (PD3O) -- node[midway, above] {$h=0$} (PDHG);
\draw[-latex] (PD3O) -- node[midway, above] {$f=0$} (LV);
\draw[-latex] (DYw)  -- node[midway, above, sloped] {completion}
  (PD3Ow);
\draw[latex-] (DYe)  -- node[midway, right]  {$A=I$} (PD3Oe);
\draw[-latex] (DRw)  -- node[midway, above, sloped] {completion} 
  (PDHGw);
\draw[latex-] (DRe)  -- node[midway, right]  {$A=I$} (PDHGe);
\draw[latex-] (PG)   -- node[midway, right]  {$A=I$} (LV);
\draw[-latex] (PD3Ow) -- node[near start, left]  {$g=0$} (PG.north);
\end{tikzpicture}
\caption{Proximal algorithms derived from PD3O.}
\label{t-sum-pd3o}
\end{figure}
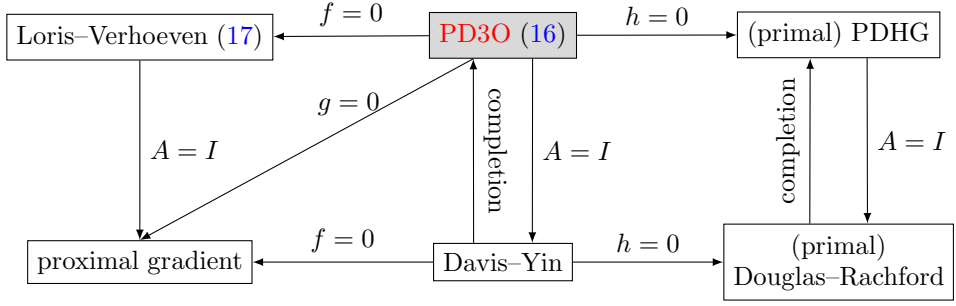

The third diagram, Figure~\ref{t-sum-pd3o}, starts with the
primal--dual three-operator (PD3O) splitting algorithm~\cite{Yan18}
\begin{subequations} \label{e-pd3o}
\begin{align}
  x^{(k+1)} &= \prox_{\tau f} (x^{(k)}-
    \tau (A^Tz^{(k)}+\nabla h(x^{(k)}))) \\
  z^{(k+1)} &= \prox_{\sigma g^\ast}
    (z^{(k)}+\sigma A(2x^{(k+1)}-x^{(k)}+\tau \nabla h(x^{(k)})
    -\tau \nabla h(x^{(k+1)}))).
\end{align}
\end{subequations}
Compared with the Condat--V\~u algorithm~\eqref{e-cv},
PD3O seems to have slightly more complicated updates and
larger complexity per iteration,
but the requirement for the stepsizes is looser:
$\sigma \tau \|A\|_2^2 \leq 1$ and $\tau \leq 1/L$.
When $h=0$, \eqref{e-pd3o} reduces to the (primal) PDHG.
The classical proximal gradient algorithm can be obtained by setting $g=0$.
When $f=0$, it reduces to the iterations%
\begin{subequations} \label{e-lv}
\begin{align}
  x^{(k+1)} &= x^{(k)}-\tau(A^Tz^{(k)}+\nabla h(x^{(k)})) \\
  z^{(k+1)} &= \prox_{\sigma g^\ast}
    \big((I-\sigma \tau AA^T)z^{(k)}
      +\sigma A(x^{(k+1)}-\tau \nabla h(x^{(k+1)})) \big).
\end{align}
\end{subequations}
This algorithm was discovered independently as
the Loris--Verhoeven algorithm~\cite{LV11},
the primal--dual fixed point algorithm
based on proximity operator (PDFP$^2$O)~\cite{CHZ13},
and the proximal alternating predictor corrector (PAPC) \cite{DST15}.
Comparison with~\eqref{e-slv} reveals a minor difference
between these two algorithms:
the gradient term in the $z$-update is taken at the newest primal
iterate $x^{(k+1)}$ in Loris--Verhoeven~\eqref{e-lv}
and at the previous point $x^{(k)}$ in the shifted version.
This difference is inherited in the proximal gradient algorithm
and its shifted version~\eqref{e-spg}.

Furthermore, when $A=I$ and $\sigma=1/\tau$ in PD3O, we recover
the well-known Davis--Yin splitting (DYS) algorithm~\cite{DY15}.
We can also set $A=I$ in~\eqref{e-lv} and obtain the iterations
\begin{subequations} \label{e-pg}
\begin{align}
  x^{(k+1)} &= x^{(k)}-\tau(z^{(k)}+\nabla h(x^{(k)})) \\
  z^{(k+1)} &= \prox_{\sigma g^\ast} \big(
    (1-\sigma\tau)z^{(k)}+\sigma(x^{(k+1)}-\tau \nabla h(x^{(k+1)})\big).
\end{align}
\end{subequations}
The stepsize conditions require $\sigma\tau \leq 1$ and $\tau \leq 1/L$.
Thus we can set $\sigma=1/\tau$ and apply Moreau decomposition.
The resulting algorithm is exactly the proximal gradient method.
The only difference in the $z$-update
between~\eqref{e-spg} and~\eqref{e-pg} is
the point at which the gradient of $h$ is taken.
The second algorithm uses the most up-to-date iterate $x^{(k+1)}$
when evaluating the gradient of $h$, and this choice allows a larger 
stepsize $\tau$.

\subsection{Primal--dual Davis--Yin (PDDY) splitting algorithm}

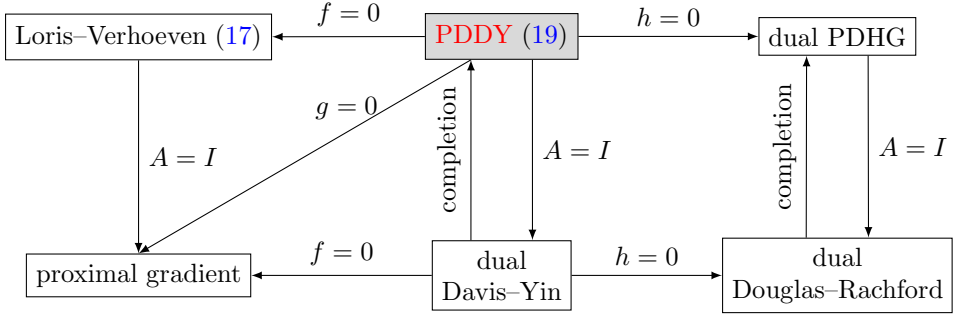
\begin{figure}
\centering
\begin{tikzpicture}[>=latex, font=\small,
                    every text node part/.style={align=center}]
\node (DY) [rectangle,draw] {dual \\ Davis--Yin};
\node (DR) [rectangle,draw, right=2.0cm of DY] 
  {dual \\ Douglas--Rachford};
\node (PG) [rectangle,draw, left=2.4cm of DY] {proximal gradient};
\node (PD3O) [rectangle,draw,above=2.4cm of DY, fill=gray!30] 
  {\red{PDDY}~\eqref{e-pddy}};
\node (PDHG) [rectangle,draw,above=2.42cm of DR] {dual PDHG};
\node (LV) [rectangle,draw,above=2.57cm of PG]
  {Loris--Verhoeven~\eqref{e-lv}};
\path (DY.north) -- (DY.north east) coordinate[pos=0.44] (DYe);
\path (DY.north) -- (DY.north west) coordinate[pos=0.44] (DYw);
\path (DR.north) -- (DR.north east) coordinate[pos=0.27] (DRe);
\path (DR.north) -- (DR.north west) coordinate[pos=0.27] (DRw);
\path (PG.north) -- (PG.north east) coordinate[pos=0.5] (PGe);
\path (PG.north) -- (PG.north west) coordinate[pos=0.5] (PGw);
\path (PD3O.south) -- (PD3O.south east) coordinate[pos=0.4] (PD3Oe);
\path (PD3O.south) -- (PD3O.south west) coordinate[pos=0.4] (PD3Ow);
\path (PDHG.south) -- (PDHG.south east) coordinate[pos=0.39] (PDHGe);
\path (PDHG.south) -- (PDHG.south west) coordinate[pos=0.39] (PDHGw);
\path (LV.south) -- (LV.south east) coordinate[pos=0.36] (LVe);
\path (LV.south) -- (LV.south west) coordinate[pos=0.36] (LVw);
\draw[-latex] (DY)   -- node[midway, above] {$h=0$} (DR);
\draw[-latex] (DY)   -- node[midway, above] {$f=0$} (PG);
\draw[-latex] (PD3O) -- node[midway, above] {$h=0$} (PDHG);
\draw[-latex] (PD3O) -- node[midway, above] {$f=0$} (LV);
\draw[-latex] (DYw)  -- node[midway, above, sloped] {completion}
  (PD3Ow);
\draw[latex-] (DYe)  -- node[midway, right]  {$A=I$} (PD3Oe);
\draw[-latex] (DRw)  -- node[midway, above, sloped] {completion} 
  (PDHGw);
\draw[latex-] (DRe)  -- node[midway, right]  {$A=I$} (PDHGe);
\draw[latex-] (PG)   -- node[midway, right]  {$A=I$} (LV);
\draw[-latex] (PD3Ow) -- node[near start, left] {$g=0$} (PG.north);
\end{tikzpicture}
\caption{Proximal algorithms derived from PDDY.}
\label{t-sum-pddy}
\end{figure}

The core algorithm in Figure~\ref{t-sum-pddy} is the primal--dual
Davis--Yin (PDDY) splitting algorithm~\cite{SCMR20}%
\begin{subequations} \label{e-pddy}
\begin{align}
  z^{(k+1)} &= \prox_{\sigma g^\ast} (z^{(k)}+\sigma Ax^{(k)}) \\
  x^{(k+1)} &= \prox_{\tau f} \big(x^{(k)}-\tau A^T(2z^{(k+1)}-z^{(k)})
    -\tau \nabla h(x^{(k)}+\tau A^T(z^{(k)}-z^{(k+1)})) \big).
\end{align}
\end{subequations}
The requirement for stepsizes is the same as that in PD3O:
$\sigma \tau \|A\|_2^2 \leq 1$ and $\tau \leq 1/L$.
Figure~\ref{t-sum-pddy} is almost identical to Figure~\ref{t-sum-pd3o}
with the roles of $f$ and $g$ exchanged.
When $h=0$, PDDY reduces to the dual PDHG.
In addition, when $A=I$ and $\sigma=1/\tau$,
PDDY reduces to the Davis--Yin algorithm, but with $f$ and $g$ exchanged.
Similarly, when $h=0$, $A=I$ and $\sigma=1/\tau$,
PDDY reverts to the Douglas--Rachford algorithm with $f$ and $g$ switched.

We have seen that the middle and right parts of Figure~\ref{t-sum-pddy}
are those of Figure~\ref{t-sum-pd3o} with $f$ and $g$ switched.
However, when one of the functions $f$ or $g$ is absent,
the algorithms reduced from PD3O and PDDY are exactly the same.
In particular, when $f=0$, PDDY reduces to the Loris--Verhoeven algorithm.

\section{Bregman distances} \label{s-bregman}

In this section we give the definition of Bregman proximal operators and 
the basic properties that will be used in the paper.  
We refer the interested reader to~\cite{CeZ:97} 
for an in-depth discussion of Bregman distances, their history, and
applications.

Let $\phi$ be a convex function with a domain that has nonempty
interior, and assume $\phi$ is continuous on $\dom \phi$ and 
continuously differentiable on $\intr(\dom \phi)$.
The \emph{generalized distance} (or \emph{Bregman distance})
generated by the \emph{kernel function} $\phi$ is defined as the function
\[
d(x,y) = \phi(x)-\phi(y)-\inprod{\nabla \phi(y)}{x-y},
\]
with domain $\dom d=\dom \phi \times \intr(\dom \phi)$.
The corresponding \emph{Bregman proximal operator} of a function $f$ is 
\begin{align}
\prox_f^\phi(y,a)
  &= \argmin_x{(f(x)+\inprod{a}{x}+d(x,y))} \label{e-prox-opt} \\
&= \argmin_x{(f(x)+\inprod{a}{x}+\phi(x)
   -\inprod{\nabla \phi(y)}{x})}. \label{e-prox-opt2}
\end{align}
It is assumed that for every $a$ and every $y \in \intr(\dom\phi)$
the minimizer $\hat x=\prox_f^\phi(y,a)$ is unique and in
$\intr(\dom \phi)$.

The distance generated by the kernel $\phi(x)=(1/2)\|x\|^2$
is the squared Euclidean distance 
$d(x,y) = (1/2) \|x-y\|^2$.  The corresponding Bregman 
proximal operator is the standard proximal operator applied to $y-a$:
\[
\prox_f^\phi(y,a)=\prox_f(y-a).
\]
For this distance, closedness and convexity of $f$ 
guarantee that the proximal operator is well defined.
The questions of existence and uniqueness are more complicated for 
general Bregman distances. 
There are no simple general conditions that guarantee that for 
every $a$ and every $y \in \intr(\dom \phi)$
the generalized proximal operator~\eqref{e-prox-opt} is uniquely
defined and in $\intr(\dom\phi)$.
Some sufficient conditions are provided
(see, for example, \cite[Section~4.1]{Bubeck15},
\cite[Assumption~A]{BBT17}),
but they may be quite restrictive or difficult to verify in practice.
In applications, however, the Bregman proximal operator is used with 
specific combinations of $f$ and $\phi$, for which the minimization problem
in~\eqref{e-prox-opt} is particularly easy to solve.
In those applications, existence and uniqueness of the solution follow
directly from the closed-form solution or availability
of a fast algorithm to compute it.
A typical example will be provided in Section~\ref{s-exp}.

From the expression~\eqref{e-prox-opt2} we see that 
$\hat x=\prox_f^\phi(y,a)$
satisfies
\[
\nabla \phi(y) - \nabla \phi(\hat x) - a \in \partial f(\hat x).
\]
Equivalently, by definition of subgradient,
\begin{align}
f(x) + \inprod{a}{x}
&\geq f(\hat x)+\inprod{a}{\hat x}
  +\inprod{\nabla\phi(y)-\nabla\phi(\hat x)}{x-\hat x} \nonumber \\
&= f(\hat x)+\inprod{a}{\hat x}+d(\hat x, y)+d(x,\hat x)-d(x,y)
  \label{e-prox-cond}
\end{align}
for all $x\in \dom f \cap \dom \phi$.

\section{Bregman Condat--V\~u three-operator splitting algorithms}
\label{s-bcv}

We now discuss two Bregman three-operator splitting algorithms
for the problem~\eqref{e-prob}.
The algorithms use a generalized distance $d_\primal$
in the primal space, generated by a kernel~$\phi_\primal$,
and a generalized distance $d_\dual$ in the dual space,
generated by a kernel $\phi_\dual$.
The first algorithm~is
\begin{subequations} \label{e-bcv}
\begin{align}
  x^{(k+1)} &= \prox^{\phi_\primal}_{\tau f}
    \big(x^{(k)}, \tau A^Tz^{(k)}+\tau \nabla h(x^{(k)}) \big) 
 \label{e-bcv-a}\\
  z^{(k+1)} &= \prox^{\phi_\dual}_{\sigma g^\ast}
    \big(z^{(k)}, -\sigma A(2x^{(k+1)}-x^{(k)})\big)
 \label{e-bcv-b}
\end{align}
\end{subequations}
and will be referred to as 
the \emph{Bregman primal Condat--V\~u algorithm}.
The second algorithm will be called
the \emph{Bregman dual Condat--V\~u algorithm}:
\begin{subequations} \label{e-bcv-dual}
\begin{align}
  z^{(k+1)} &= \prox_{\sigma g^\ast}^{\phi_\dual}
    (z^{(k)}, -\sigma Ax^{(k)}) 
 \label{e-bcv-dual-a}\\
  x^{(k+1)} &= \prox_{\tau f}^{\phi_\primal}
    (x^{(k)}, \tau A^T(2z^{(k+1)}-z^{(k)}) + \tau\nabla h(x^{(k)}) ).
 \label{e-bcv-dual-b}
\end{align}
\end{subequations}
The two algorithms need starting points
$x^{(0)} \in \intr(\dom\phi_\primal) \cap \dom h$,
and $z^{(0)} \in \intr(\dom\phi_\dual)$.
Conditions on stepsizes $\sigma$, $\tau$ will be specified later.
When Euclidean distances are used for the primal and dual proximal
operators, the two algorithms reduce to the primal and dual variants
of the Condat--V\~u algorithm~\eqref{e-cv} and~\eqref{e-cv-dual},
respectively.
Algorithm~\eqref{e-bcv} has been proposed in~\cite{ChP:16}.  Here
we discuss it together with~\eqref{e-bcv-dual} in a unified framework.

In Section~\ref{s-bcv-ppm} we show that the proposed algorithms
can be interpreted as the Bregman proximal point method
applied to a monotone inclusion problem.
In Section~\ref{s-bcv-conv} we analyze their convergence.
In Section~\ref{s-bcv-connection} we discuss the connections between the 
two algorithms and other Bregman proximal splitting methods.

\paragraph{Assumptions} 
Throughout Section~\ref{s-bcv} we make the following assumptions.
The kernel functions $\phi_\primal$ and $\phi_\dual$ are 
1-strongly convex with respect to norms $\|\cdot\|_\primal$ and 
$\|\cdot\|_\dual$, respectively:
\begin{equation} \label{e-bcv-str-cvx}
d_\primal(x,x^\prime) \geq \frac{1}{2}\|x-x^\prime\|^2_\primal, \qquad
d_\dual(z,z^\prime) \geq \frac{1}{2}\|z-z^\prime\|^2_\dual
\end{equation}
for all $(x,x') \in \dom d_\primal$ and $(z,z') \in \dom d_\dual$.
The assumption that the strong convexity constants are equal to one
can be made without loss of generality, by scaling the norms (or distances)
if needed.
We also assume that the function $L\phi_\primal-h$ is convex for some 
$L>0$.  More precisely, $\dom \phi_\primal \subseteq \dom h$ and
\BEQ \label{e-bcv-lip}
h(x)-h(x^\prime)-\inprod{\nabla h(x^\prime)}{x-x^\prime}
\leq Ld_\primal (x,x^\prime) \quad \mbox{for all\ }
(x,x^\prime) \in \dom d_\primal.
\EEQ
Note that this assumption is looser than
the one in~\cite[Equation~(4)]{ChP:16}.
We denote by $\|A\|$ the matrix norm
\BEQ \label{e-norm}
\|A\|=\sup_{u \neq 0, v \neq 0}
    \frac{\inprod{v}{Au}}{\|v\|_\dual \|u\|_\primal}
 = \sup_{u\neq 0} \frac{\|Au\|_{\dual,*}}{\|u\|_\primal}
 = \sup_{v\neq 0} \frac{\|A^Tv\|_{\primal,*}}{\|v\|_\dual},
\EEQ
where $\|\cdot\|_{\primal,*}$ and $\|\cdot\|_{\dual,*}$ are
the dual norms of $\|\cdot\|_\primal$ and $\|\cdot\|_\dual$. 

It is also assumed that the primal--dual optimality 
conditions~\eqref{e-opt-cond} have a solution $(x^\star, z^\star)$
with $x^\star \in \dom \phi_\primal$ and 
$z^\star \in \dom \phi_\dual$.  

\subsection{Derivation from Bregman proximal point method}
\label{s-bcv-ppm}

The Bregman Condat--V\~u algorithms~\eqref{e-bcv} and~\eqref{e-bcv-dual} 
can be viewed as applications of the Bregman proximal point algorithm
to the optimality conditions~\eqref{e-opt-cond}.
This interpretation extends the derivation of the Bregman PDHG algorithm 
from the Bregman proximal point algorithm given in \cite{JV22}.
The idea originates with He and Yuan's interpretation of PDHG as a 
``preconditioned'' proximal point algorithm~\cite{HeY:12}.

The Bregman proximal point algorithm~\cite{Eck:93,CeZ:97,Gul:94}
is an algorithm for monotone inclusion problems $0 \in F(u)$.  
The update $u^{(k+1)}$ in one iteration of the algorithm is defined 
as the solution of the inclusion
\[
\nabla \phi(u^{(k)}) - \nabla \phi(u^{(k+1)}) \in F(u^{(k+1)}),
\]
where $\phi$ is a Bregman kernel function.
Applied to~\eqref{e-opt-cond}, with a kernel function
$\phi_\mathrm{pd}$, the algorithm generates a sequence 
$(x^{(k)}, z^{(k)})$ defined by 
\begin{equation} \label{e-bcv-ppa}
\nabla \phi_\mathrm{pd}(x^{(k)},z^{(k)}) 
-\nabla \phi_\mathrm{pd}(x^{(k+1)},z^{(k+1)})
\in 
\begin{bmatrix} A^T z^{(k+1)} +  
\partial f(x^{(k+1)})+\nabla h(x^{(k+1)}) \\ 
-A x^{(k+1)} + \partial g^\ast(z^{(k+1)}) \end{bmatrix}.
\end{equation}

\subsubsection{Primal--dual Bregman distances} \label{s-pd-dist}

We introduce four possible primal--dual kernel functions: the functions
\[
\phi_+(x,z) = \frac{1}{\tau}\phi_\primal(x)
+\frac{1}{\sigma}\phi_\dual(z)+\inprod{z}{Ax}, \qquad
\phi_-(x,z) = \frac{1}{\tau}\phi_\primal(x)
+\frac{1}{\sigma}\phi_\dual(z)-\inprod{z}{Ax},
\]
where $\sigma,\tau > 0$, and the functions 
\[
 \phi_\mathrm{dcv}(x,z) = \phi_+(x,z) - h(x), \qquad 
 \phi_\mathrm{pcv}(x,z) = \phi_-(x,z) - h(x).
\]
The subscripts in $\phi_+$ and $\phi_-$ 
refer to the sign of the inner product term $\inprod{z}{Ax}$.
The subscripts in $\phi_\mathrm{pcv}$ and $\phi_\mathrm{dcv}$
indicate the algorithm (Bregman primal or dual Condat-V\~u) for which
these distances will be relevant. 
If these kernel functions are convex,
they generate the following Bregman distances.
The distances generated by $\phi_+$ and $\phi_-$ are
\begin{align}
d_+ (x,z;x^\prime,z^\prime) &= 
  \frac{1}{\tau}d_\primal(x,x^\prime)+\frac{1}{\sigma}d_\dual(z,z^\prime)
  +\inprod{z-z^\prime}{A(x-x^\prime)} \nonumber \\
d_-(x,z;x^\prime,z^\prime) &=
  \frac{1}{\tau}d_\primal(x,x^\prime)+\frac{1}{\sigma}d_\dual(z,z^\prime)
  -\inprod{z-z^\prime}{A(x-x^\prime)}, \label{e-bcv-dminus}
\end{align}
respectively, and the distances generated by $\phi_\mathrm{dcv}$ and 
$\phi_\mathrm{pcv}$ are  
\begin{align*}
d_\mathrm{dcv}(x,z;x^\prime,z^\prime) &= 
  d_+ (x,z;x^\prime,z^\prime) 
  -h(x)+h(x^\prime)+\inprod{\nabla h(x')}{x-x^\prime} \\
d_\mathrm{pcv}(x,z;x^\prime,z^\prime) &= 
  d_- (x,z;x^\prime,z^\prime) 
  -h(x)+h(x^\prime)+\inprod{\nabla h(x')}{x-x^\prime}.
\end{align*}
We now show that $\phi_+$ and $\phi_-$ are convex if
\[
\sigma\tau \|A\|^2 \leq 1 
\]
and strongly convex if $\sigma\tau\|A\|^2 < 1$,
and that the functions $\phi_\mathrm{dcv}$ and 
$\phi_\mathrm{pcv}$ are convex if
\BEQ \label{e-bcv-param}
\sigma\tau \|A\|^2 + \tau L \leq 1
\EEQ
and strongly convex if $\sigma\tau \|A\|^2 + \tau L < 1$.
\begin{proof}
To show that the kernel functions $\phi_+$ and $\phi_-$ are convex,
we show that $d_+$ and $d_-$ are nonnegative.
Suppose $\sigma\tau \|A\|^2 \leq \delta_1 \delta_2$
with $\delta_1,\delta_2 > 0$.
Then~\eqref{e-bcv-str-cvx} and the arithmetic--geometric mean inequality
imply that
\begin{align}
\left| \inprod{z-z^\prime}{A(x-x^\prime)}\right|
  &\leq \|A\| \|z-z'\|_\dual \|x-x'\|_\primal  \nonumber \\
&\leq \sqrt{\frac{\delta_1 \delta_2}{\sigma \tau}}
\|z-z'\|_\dual \|x-x'\|_\primal \nonumber \\
&\leq \frac{\delta_1}{2\tau} \|x-x'\|_\primal^2  +
  \frac{\delta_2}{2\sigma} \|z-z'\|_\primal^2 \nonumber \\
&\leq \frac{\delta_1}{\tau} d_\primal (x,x') +
  \frac{\delta_2}{\sigma} d_\dual (z,z'). \label{e-bcv-delta}
\end{align}
Therefore,
\begin{align*}
d_{\pm} (x,z; x',z')
  &= \frac{1}{\tau} d_\primal(x,x^\prime) + \frac{1}{\sigma}
  d_\dual(z,z^\prime) \pm \inprod{z-z^\prime}{A(x-x^\prime)} \\
&\geq \frac{1-\delta_1}{\tau} d_\primal (x,x') +
  \frac{1-\delta_2}{\sigma} d_\dual (z,z') \\
&\geq \frac{1-\delta_1}{2\tau} \|x-x'\|_\primal^2 +
  \frac{1-\delta_2}{2\sigma} \|z-z'\|_\dual^2.
\end{align*}
With $\delta_1=\delta_2=1$, this shows convexity of $\phi_+$ and $\phi_-$;
with $\delta_1=\delta_2<1$, strong convexity.
Similarly,
\begin{align*}
d_\mathrm{dcv/pcv} (x,z; x', z')
  &= d_{\pm} (x,z; x',z') - h(x) + h(x') + \inprod{\nabla h(x')}{x-x'} \\
&\geq \frac{1-\tau L-\delta_1}{\tau} d_\primal (x,x^\prime) +
  \frac{1-\delta_2}{\sigma} d_\dual (z,z^\prime).
\end{align*}
With $\delta_1=1-\tau L$ and $\delta_2=1$, this shows convexity of
$\phi_\mathrm{pcv}$ and $\phi_\mathrm{dcv}$;
with $\delta_1=\delta-\tau L$ and $\delta_2=\delta < 1$,
strong convexity.
\end{proof}

\subsubsection{Bregman Condat-V\~u algorithms from proximal point method}

The Bregman primal Condat--V\~u algorithm~\eqref{e-bcv}
is the Bregman proximal point method with the kernel function 
$\phi_\mathrm{pd} = \phi_\mathrm{pcv}$.
If we take $\phi_\mathrm{pd} = \phi_\mathrm{pcv}$ in~\eqref{e-bcv-ppa},
we obtain two coupled inclusions that determine $x^{(k+1)}$, $z^{(k+1)}$.  
The first one is
\BEAS
0 &\in & \frac{1}{\tau}(\nabla \phi_\primal(x^{(k+1)})
  -\nabla \phi_\primal(x^{(k)}))-A^T(z^{(k+1)}-z^{(k)})
  -\nabla h(x^{(k+1)})+\nabla h(x^{(k)}) \\
&  & \mbox{} +A^Tz^{(k+1)}+\partial f(x^{(k+1)})+\nabla h(x^{(k+1)}) \\
&= & \frac{1}{\tau}(\nabla \phi_\primal(x^{(k+1)})
  -\nabla \phi_\primal(x^{(k)}))+A^Tz^{(k)}
  +\nabla h(x^{(k)})+\partial f(x^{(k+1)}).
\EEAS
This shows that $x^{(k+1)}$ solves the optimization problem
\[
  \mbox{minimize} \quad f(x)+\inprod{A^Tz^{(k)}+\nabla h(x^{(k)})}{x}
  +\frac{1}{\tau} d_\primal(x,x^{(k)}).
\]
The solution is the $x$-update~\eqref{e-bcv-a} in the Bregman primal
Condat--V\~{u} method. The second inclusion~is
\BEAS
0 &\in & \frac{1}{\sigma}(\nabla \phi_\dual(z^{(k+1)})
  -\nabla \phi_\dual(z^{(k)}))-A(x^{(k+1)}-x^{(k)})-Ax^{(k+1)}
  +\partial g^\ast(z^{(k+1)}) \\
&= & \frac{1}{\sigma}(\nabla \phi_\dual(z^{(k+1)})
  -\nabla \phi_\dual(z^{(k)}))-A(2x^{(k+1)}-x^{(k)})
  +\partial g^\ast(z^{(k+1)}).
\EEAS
This shows that $z^{(k+1)}$ solves the optimization problem
\[
  \mbox{minimize} \quad g^\ast(z)-\inprod{z}{A(2x^{(k+1)}-x^{(k)})}
  +\frac{1}{\sigma} d_\dual(z,z^{(k)}).
\]
The solution is the $z$-update~\eqref{e-bcv-b}.

Choosing $\phi_\mathrm{pd} = \phi_\mathrm{dcv}$ in~\eqref{e-bcv-ppa}
yields the Bregman dual Condat--V\~u algorithm~\eqref{e-bcv-dual}.
Substituting $\phi_\mathrm{pd} = \phi_\mathrm{dcv}$ in~\eqref{e-bcv-ppa}
gives the inclusions
\BEAS
0 &\in & \frac{1}{\tau}(\nabla \phi_\primal(x^{(k+1)})
  -\nabla \phi_\primal(x^{(k)}))+A^T(z^{(k+1)}-z^{(k)})
  -\nabla h(x^{(k+1)})+\nabla h(x^{(k)}) \\
&  & \mbox{} +A^Tz^{(k+1)}+\partial f(x^{(k+1)})+\nabla h(x^{(k+1)}) \\
&= & \frac{1}{\tau}(\nabla \phi_\primal(x^{(k+1)})
  -\nabla \phi_\primal(x^{(k)}))+A^T(2z^{(k+1)} - z^{(k)})
  +\nabla h(x^{(k)})+\partial f(x^{(k+1)})
\EEAS
and
\BEAS
0 &\in & \frac{1}{\sigma}(\nabla \phi_\dual(z^{(k+1)})
  -\nabla \phi_\dual(z^{(k)}))+A(x^{(k+1)}-x^{(k)})-Ax^{(k+1)}
  +\partial g^\ast(z^{(k+1)}) \\
&= & \frac{1}{\sigma}(\nabla \phi_\dual(z^{(k+1)})
  -\nabla \phi_\dual(z^{(k)}))-A x^{(k)} +\partial g^\ast(z^{(k+1)}).
\EEAS
The second inclusion shows that $z^{(k+1)}$ is given by the
$z$-update~\eqref{e-bcv-dual-a}.
Given $z^{(k+1)}$, one can solve the first inclusion for $x^{(k+1)}$
and obtains the $x$-update~\eqref{e-bcv-dual-b}.

\subsection{Convergence analysis} \label{s-bcv-conv}

The derivation in Section~\ref{s-bcv-ppm} allows us to apply
existing convergence theory for the Bregman proximal point method
to the proposed algorithms~\eqref{e-bcv} and~\eqref{e-bcv-dual}.
In particular, Solodov and Svaiter~\cite{SS00} have studied
Bregman proximal point methods with inexact prox-evaluations
for solving variational inequalities,
which include the monotone inclusion problem as a special case.
The results in~\cite{SS00} can be applied to analyze convergence
of the Bregman Condat--V\~u methods with inexact evaluations of
proximal operators.

The literature on the Bregman proximal point method
for monotone inclusions~\cite{Eck:93,Gul:94,SS00}
focuses on the convergence of iterates,
and this generally requires additional assumptions on 
$\phi_\primal$ and $\phi_\dual$
(beyond the assumptions of convexity made in~Section~\ref{s-bcv-ppm}).
In this section we present a self-contained convergence analysis
and give a direct proof of an $O(1/k)$ rate of ergodic convergence.
We also give a self-contained proof of convergence of the 
iterates $x^{(k)}$ and~$z^{(k)}$.

We make the assumptions listed in~Section~\ref{s-bcv-ppm}:
the strong convexity assumption~\eqref{e-bcv-str-cvx} for the primal
and dual kernels $\phi_\primal$ and $\phi_\dual$, 
and the \emph{relative smoothness} property~\eqref{e-bcv-lip}
of the function $h$.  We assume that the stepsizes $\sigma$, $\tau$
satisfy~\eqref{e-bcv-param}, and that the primal--dual optimality
condition~\eqref{e-opt-cond} has a solution 
$(x^\star, z^\star) \in \dom \phi_\primal \times \dom \phi_\dual$.

For the sake of brevity we combine the analysis of the Bregman primal 
and the Bregman dual Condat-V\~u algorithms.
In the following, $d$, $\tilde d$, $\tilde \phi$ are defined as
\[
\begin{array}{llll}
d = d_- \quad & \tilde d = d_\mathrm{pcv} \quad
  & \tilde \phi = \phi_\mathrm{pcv} \qquad
  & \mbox{for Bregman primal Condat--V\~u~\eqref{e-bcv},} \\
d = d_+ \quad &\tilde d = d_\mathrm{dcv} \quad
  & \tilde \phi = \phi_\mathrm{dcv} \qquad
  & \mbox{for Bregman dual Condat--V\~u~\eqref{e-bcv-dual}.}
\end{array}
\]

\subsubsection{One-iteration analysis}

We first show that the iterates 
$x^{(k+1)}$, $z^{(k+1)}$ generated by the Bregman Condat--V\~u
algorithms~\eqref{e-bcv} and~\eqref{e-bcv-dual} satisfy
\BEA
\lefteqn{ \cL(x^{(k+1)},z)-\cL(x,z^{(k+1)})} \nonumber \\
& \leq & d(x,z;x^{(k)},z^{(k)})- d(x,z;x^{(k+1)},z^{(k+1)}) 
 - \tilde d(x^{(k+1)}, z^{(k+1)}; x^{(k)}, z^{(k)})
 \label{e-bcv-conv-i} 
\EEA
for all $x \in \dom f \cap \dom \phi_\primal$ and 
$z \in \dom g^\ast \cap \dom \phi_\dual$.  

\begin{proof}
We write~\eqref{e-bcv} and~\eqref{e-bcv-dual} in a unified notation as
\begin{subequations}
\begin{align}
x^{(k+1)} &= \prox_{\tau f}^{\phi_\primal} (x^{(k)},
  \tau (A^T \tilde z + \nabla h(x^{(k)}))) \label{e-bcv-conv-i-primal} \\
z^{(k+1)} &= \prox_{\sigma g^*}^{\phi_\dual}
  (z^{(k)}, -\sigma A\tilde x) \label{e-bcv-conv-i-dual}
\end{align}
\end{subequations}
where $\tilde x$ and $\tilde z$ are defined in the following table:
\[
\begin{array}{lll}
\mbox{Bregman primal Condat--V\~{u} algorithm} & 
  \tilde x=2x^{(k+1)}-x^{(k)} & \tilde z=z^{(k)} \\
\mbox{Bregman dual Condat--V\~{u} algorithm} & \tilde x=x^{(k)}
  & \tilde z=2z^{(k+1)}-z^{(k)}.
\end{array}
\]
The optimality condition~\eqref{e-prox-cond} for the proximal
operator evaluation~(\ref{e-bcv-conv-i-primal}) is that
\[
\tau (f(x^{(k+1)})- f(x)) \leq
d_\primal(x,x^{(k)})-d_\primal(x^{(k+1)},x^{(k)})-d_\primal(x,x^{(k+1)}) 
  + \tau \inprod{A^T\tilde z+\nabla h(x^{(k)})}{x-x^{(k+1)}}
\]
for all $x \in \dom f \cap \dom \phi_\primal$.
The optimality condition for~\eqref{e-bcv-conv-i-dual} is that
\[
\sigma (g^\ast(z^{(k+1)})- g^\ast(z)) \leq
  d_\dual(z,z^{(k)}) - d_\dual(z^{(k+1)}, z^{(k)}) - d_\dual(z, z^{(k+1)})
  - \sigma \inprod{z-z^{(k+1)}}{A\tilde x}
\]
for all $z \in \dom g^\ast \cap \dom \phi_\dual$.
Combining the two inequalities gives
\BEA
\lefteqn{\cL(x^{(k+1)}, z) - \cL(x,z^{(k+1)})} \nonumber \\
& = & f(x^{(k+1)}) - f(x) + h(x^{(k+1)})  - h(x) 
   + g^*(z^{(k+1)}) - g^*(z) + \inprod{A^Tz}{x^{(k+1)}}
  - \inprod{z^{(k+1)}}{Ax} \nonumber \\
&\leq &\mbox{} \frac{1}{\tau} \Big(
  d_\primal(x,x^{(k)})-d_\primal(x,x^{(k+1)})
  - d_\primal(x^{(k+1)}, x^{(k)}) \Big) \nonumber \\
& &\mbox{} + \frac{1}{\sigma} \Big(d_\dual(z,z^{(k)})-d_\dual(z,z^{(k+1)})
  - d_\dual(z^{(k+1)}, z^{(k)}) \Big) \nonumber \\
& & \mbox{} + h(x^{(k+1)})-h(x)+\inprod{\nabla h(x^{(k)})}{x-x^{(k+1)}}
  \nonumber \\
& &\mbox{} 
  + \inprod{A^T \tilde z}{x-x^{(k+1)}}
  - \inprod{z-z^{(k+1)}}{A\tilde x}
  + \inprod{A^T z}{x^{(k+1)}}-\inprod{z^{(k+1)}}{Ax}
   \label{e-bcv-conv-prf-2} \\
&\leq &\mbox{} \frac{1}{\tau} \Big(
  d_\primal(x,x^{(k)})-d_\primal(x,x^{(k+1)})
  - d_\primal(x^{(k+1)}, x^{(k)}) \Big)  \nonumber \\
& &\mbox{} + \frac{1}{\sigma} \Big(d_\dual(z,z^{(k)})-d_\dual(z,z^{(k+1)})
  - d_\dual(z^{(k+1)}, z^{(k)}) \Big) \nonumber \\
& & \mbox{} + h(x^{(k+1)})-h(x^{(k)})-
 \inprod{\nabla h(x^{(k)})}{x^{(k+1)} - x^{(k)}} \nonumber \\
& &\mbox{} + \inprod{A^T \tilde z}{x-x^{(k+1)}}
  - \inprod{z-z^{(k+1)}}{A\tilde x}
  + \inprod{A^Tz}{x^{(k+1)}}-\inprod{z^{(k+1)}}{Ax}
 \label{e-bcv-conv-prf-1}
\EEA
for all $x\in \dom f \cap\dom\phi_\primal$ and all
$z\in \dom g^* \cap \dom \phi_\dual$.
The second inequality follows from convexity of $h$.
Substituting the expressions for $\tilde x$ and $\tilde z$ in 
the Bregman primal Condat--V\~u algorithm~\eqref{e-bcv}, 
we obtain for the last line of~\eqref{e-bcv-conv-prf-1}
\BEAS
\lefteqn{\inprod{A^T\tilde z}{x-x^{(k+1)}}
  -\inprod{z-z^{(k+1)}}{A\tilde x} 
  +\inprod{A^T z}{x^{(k+1)}} - \inprod{z^{(k+1)}}{Ax}} \nonumber \\
&=& \inprod{z^{(k)}}{A(x-x^{(k+1)})}
  -\inprod{z-z^{(k+1)}}{A(2x^{(k+1)}-x^{(k)})} 
  +\inprod{A^T z}{x^{(k+1)}} - \inprod{z^{(k+1)}}{Ax} \nonumber \\
&=& \inprod{z^{(k)}-z^{(k+1)}}{A(x-x^{(k+1)})}
  +\inprod{z-z^{(k+1)}}{A(x^{(k)}-x^{(k+1)})} \nonumber \\
&=& -\inprod{z-z^{(k)}}{A(x-x^{(k)})}
  +\inprod{z-z^{(k+1)}}{A(x-x^{(k+1)})} 
+\inprod{z^{(k+1)}-z^{(k)}}{A(x^{(k+1)}-x^{(k)})}.
\EEAS
If we substitute the expressions for $\tilde x$ and $\tilde z$
in the Bregman dual Condat--V\~u algorithm, the last line
of~(\ref{e-bcv-conv-prf-1}) becomes
\BEAS
\lefteqn{\inprod{A^T\tilde z}{x-x^{(k+1)}}
  -\inprod{z-z^{(k+1)}}{A\tilde x}
  +\inprod{A^T z}{x^{(k+1)}} - \inprod{z^{(k+1)}}{Ax}} \\
&=& \inprod {A^T(z-z^{(k)})} {x-x^{(k)}}
  -\inprod {A^T(z-z^{(k+1)})} {x-x^{(k+1)}}
  -\inprod {A^T(z^{(k+1)}-z^{(k)})} {x^{(k+1)}-x^{(k)}}.
\EEAS
Therefore, for both algorithms, \eqref{e-bcv-conv-prf-1}  implies that
\BEAS
\lefteqn{\cL(x^{(k+1)}, z) - \cL(x,z^{(k+1)})} \nonumber \\
&\leq& \frac{1}{\tau} d_\primal(x,x^{(k)})
  + \frac{1}{\sigma} d_\dual(z,z^{(k)})
  \mp \inprod{z-z^{(k)}}{A(x-x^{(k)})} \nonumber \\
& &\mbox{} - \Big(\frac{1}{\tau} d_\primal(x,x^{(k+1)})
  + \frac{1}{\sigma} d_\dual(z,z^{(k+1)})
  \mp \inprod{z-z^{(k+1)}}{A(x-x^{(k+1)})} \Big) \nonumber \\
& &\mbox{} -\Big(\frac{1}{\tau} d_\primal(x^{(k+1)},x^{(k)})
  + \frac{1}{\sigma} d_\dual(z^{(k+1)},z^{(k)})
  \mp \inprod{z^{(k+1)} - z^{(k)}}{A(x^{(k+1)}-x^{(k)})} \Big) 
\nonumber \\
& & \mbox{} +
  h(x^{(k+1)})-h(x^{(k)}) -\inprod{\nabla h(x^{(k)})}{x^{(k+1)}-x^{(k)}},
\EEAS
if we select the minus sign in $\mp$ for the Bregman primal Condat--V\~u
algorithm, and the plus sign for the Bregman dual Condat--V\~u algorithm.
For the primal method, this shows
\BEAS
\lefteqn{\cL(x^{(k+1)}, z) - \cL(x,z^{(k+1)})} \\
 &\leq & d_-(x,z;x^{(k)},z^{(k)})- d_-(x,z;x^{(k+1)},z^{(k+1)}) 
 - d_\mathrm{pcv}(x^{(k+1)},z^{(k+1)};x^{(k)},z^{(k)}).
\EEAS
For the dual method,
\BEAS
\lefteqn{\cL(x^{(k+1)}, z) - \cL(x,z^{(k+1)})} \\
&\leq & d_+(x,z;x^{(k)},z^{(k)})- d_+(x,z;x^{(k+1)},z^{(k+1)}) 
  - d_\mathrm{dcv}(x^{(k+1)},z^{(k+1)};x^{(k)},z^{(k)}). 
\EEAS
\end{proof}

\subsubsection{Ergodic convergence}

We define averaged iterates 
\begin{equation} \label{e-bcv-avg}
x^{(k)}_\avg = \frac{1}{k} \sum_{i=1}^k x^{(i)}, \qquad 
z^{(k)}_\avg = \frac{1}{k} \sum_{i=1}^k z^{(i)} 
\end{equation} 
for $k \geq 1$.
We show that 
\begin{equation} \label{e-bcv-conv}
\cL(x^{(k)}_\avg, z) - \cL(x, z^{(k)}_\avg)
\leq \frac{2}{k} \Big(\frac{1}{\tau} d_\primal (x,x^{(0)})
  + \frac{1}{\sigma} d_\dual (z,z^{(0)}) \Big)
\end{equation}
for all 
$x \in \dom f \cap \dom \phi_\primal$ and
$z \in \dom g^\ast \cap \dom \phi_\dual$.

\begin{proof}
From~\eqref{e-bcv-conv-i},
since $\cL(u,v)$ is convex in $u$ and concave in $v$,
\begin{align*}
\cL(x^{(k)}_\avg, z) - \cL(x, z^{(k)}_\avg)
  &\leq \frac{1}{k} \sum_{i=1}^k
  \big(\cL(x^{(i)},z)-\cL(x,z^{(i)})\big) \nonumber \\
&\leq \frac{1}{k} \big(d(x,z;x^{(0)},z^{(0)})
  - d(x,z;x^{(k)},z^{(k)}) \big) \nonumber \\
&\leq \frac{1}{k} d(x,z;x^{(0)},z^{(0)}) \nonumber \\
&\leq \frac{2}{k} \Big(\frac{1}{\tau} d_\primal (x,x^{(0)})
  + \frac{1}{\sigma} d_\dual (z,z^{(0)}) \Big)
\end{align*}
for all $x \in \dom f \cap \dom \phi_\primal$
and $z \in \dom g^\ast \cap \dom \phi_\dual$.
The last step follows from~\eqref{e-bcv-delta} with $\delta_1=\delta_2=1$.
\end{proof}

Substituting $x=x^\star$, $z=z^\star$ in~\eqref{e-bcv-conv} gives
\[
\cL(x^{(k)}_\avg, z^\star)-\cL(x^\star, z^{(k)}_\avg)
\leq \frac{2}{k} \Big(\frac{1}{\tau} d_\primal (x^\star,x^{(0)})
  + \frac{1}{\sigma} d_\dual (z^\star,z^{(0)}) \Big).
\]
More generally, if $X \subseteq \dom \phi_\primal$ and
$Z \subseteq \dom \phi_\dual$ are compact convex sets
that contain optimal solutions $x^\star$, $z^\star$ in their interiors,
then the merit function~\eqref{e-merit} is bounded by
\[
\eta(x^{(k)}_\avg, z^{(k)}_\avg)
 \leq \frac{2}{k} 
  \Big(\frac{1}{\tau} \sup_{x \in X} d_\primal(x,x^{(0)})
   + \frac{1}{\sigma} \sup_{z \in Z} d_\dual(z,z^{(0)}) \Big).
\]

\subsubsection{Monotonicity properties} \label{s-fejer}

For $x=x^\star$, $z= z^\star$, 
the left-hand side of~\eqref{e-bcv-conv-i} is nonnegative and therefore
\BEQ
d(x^\star, z^\star; x^{(k+1)}, z^{(k+1)}) 
\leq d(x^\star, z^\star; x^{(k)}, z^{(k)})  
  -\tilde d(x^{(k+1)}, z^{(k+1)}; x^{(k)}, z^{(k)}) 
  \label{e-bcv-fejer}
\EEQ
for $k \geq 0$. Hence
$d(x^\star, z^\star; x^{(k+1)}, z^{(k+1)}) \leq
d(x^\star, z^\star; x^{(k)}, z^{(k)})$ and
\BEQ
d(x^\star, z^\star; x^{(k)}, z^{(k)})
\leq 
d(x^\star, z^\star; x^{(0)}, z^{(0)}). \label{e-bcv-bounded}
\EEQ
The inequality~\eqref{e-bcv-fejer} also implies that
\[
\sum_{i=0}^k \tilde d(x^{(i+1)}, z^{(i+1)}; x^{(i)}, z^{(i)})
\leq d(x^\star, z^\star; x^{(0)}, z^{(0)}).
\] 
Hence $\tilde d(x^{(k+1)}, z^{(k+1)}; x^{(k)}, z^{(k)}) \rightarrow 0$.

\subsubsection{Convergence of iterates}

Convergence of iterates can be obtained by combining the derivation
in Section~\ref{s-bcv-ppm} and existing results on Bregman proximal point
method \cite[Theorem~3.1]{Gul:94}, \cite[Theorem~3.2]{SS00}.
Here we provide a self-contained proof under additional assumptions 
about the primal and dual distance functions.
The following two assumptions are common in the literature
on Bregman distances~\cite{ChT:93,Eck:93,Gul:94,CeZ:97}.
\begin{enumerate}
\item For fixed $x$ and $z$, the sublevel sets
$\{x^\prime \mid d_\primal(x,x^\prime) \leq \gamma\}$ and
$\{z^\prime \mid d_\dual(z,z^\prime) \leq \gamma\}$
are closed.
In other words, the distances $d_\primal(x,x^\prime)$ 
and $d_\dual(z,z^\prime)$ are closed functions of 
$x^\prime$ and $z^\prime$, respectively.
Since a sum of closed functions is closed, the distance
$d(x,z; x',z')$ is a closed function of $(x',z')$, for fixed $(x,z)$.

\item If $\tilde x^{(k)} \in \intr(\dom \phi_\primal)$ converges to
$x \in \dom \phi_\primal$, then $d_\primal(x,\tilde x^{(k)}) \to 0$.
Similarly, if $\tilde z^{(k)} \in \intr(\dom \phi_\dual)$ converges to
$z \in \dom \phi_\dual$, then $d_\dual(z,\tilde z^{(k)}) \to 0$.
\end{enumerate}
We also assume that 
$\sigma\tau \|A\|^2 + \tau L < 1$.
As shown in Section~\ref{s-pd-dist} this implies that
the kernel functions $\phi_\mathrm{pcv}$ and $\phi_\mathrm{dcv}$ are
strongly convex and that
\BEQ \label{e-bcv-str-cvx-1}
 \tilde d(x,z; x',z') \geq 
 \frac{\alpha}{2\tau} \|x-x'\|_\primal^2
 + \frac{\alpha}{2\sigma} \|z-z'\|_\dual^2
\EEQ
for some $\alpha > 0$.  Similarly, $\sigma\tau\|A\|^2 < 1$ implies 
that
\BEQ \label{e-bcv-str-cvx-2}
 d(x,z; x',z') \geq 
 \frac{\beta}{2\tau} \|x-x'\|_\primal^2
 + \frac{\beta}{2\sigma} \|z-z'\|_\dual^2
\EEQ
for some $\beta > 0$.
Recall that $d=d_-$, $\tilde d=d_\mathrm{pcv}$
for the Bregman primal Condat--V\~u algorithm~\eqref{e-bcv},
and $d=d_+$, $\tilde d=d_\mathrm{dcv}$ 
for the Bregman dual Condat--V\~u
algorithm.
\begin{proof}
We first note that  
$\tilde d(x^{(k+1)}, z^{(k+1)}; x^{(k)}, z^{(k)}) \rightarrow 0$
and~(\ref{e-bcv-str-cvx-1})
imply that $x^{(k+1)} -x^{(k)} \rightarrow 0$ and
$z^{(k+1)} -z^{(k)} \rightarrow 0$.

The inequality~\eqref{e-bcv-bounded}, together with~\eqref{e-bcv-str-cvx-2},
implies that the sequence $(x^{(k)},z^{(k)})$ is bounded.
Let $(x^{(k_i)},z^{(k_i)})$ be a convergent subsequence of 
$(x^{(k)}, z^{(k)})$ with limit $(\hat x, \hat z)$.
Since $x^{(k_i+1)}-x^{(k_i)} \rightarrow 0$ and
$z^{(k_i+1)} -z^{(k_i)} \rightarrow 0$, the sequence
$(x^{(k_i+1)},z^{(k_i+1)})$ also converges to $(\hat x, \hat z)$.
We show that $(\hat x,\hat z)$ satisfies the optimality condition
\eqref{e-opt-cond}.

From~\eqref{e-bcv-bounded}, $d(x^\star,z^\star;x^{(k_i)},z^{(k_i)})$
is bounded.
Since the sublevel sets $\{(x^\prime,z^\prime) \mid d(x^\star,z^\star;
x^\prime,z^\prime) \leq \gamma\}$ are closed subsets of
$\intr(\dom \phi_\primal) \cap \intr(\dom \phi_\dual)$,
the limit $(\hat x,\hat z) \in
\intr(\dom \phi_\primal) \cap \intr(\dom \phi_\dual)$.
The iterates in the subsequence satisfy
\begin{equation} \label{e-bcv-fejer-prf}
\nabla \phi_\mathrm{pd} (x^{(k_i)},z^{(k_i)})
-\nabla \phi_\mathrm{pd} (x^{(k_i+1)},z^{(k_i+1)})
+\begin{bmatrix} -A^Tz^{(k_i+1)} \\ Ax^{(k_i+1)} \end{bmatrix} \in
\begin{bmatrix} \partial f(x^{(k_i+1)})+\nabla h(x^{(k_i+1)}) \\
\partial g^\ast(z^{(k_i+1)}) \end{bmatrix},
\end{equation} 
where $\phi_\mathrm{pd}=\phi_\mathrm{pcv}$ in the Bregman primal
Condat--V\~u algorithm and $\phi_\mathrm{pd}=\phi_\mathrm{dcv}$
in the Bregman dual Condat--V\~u algorithm.
The left-hand side of~\eqref{e-bcv-fejer-prf}
converges to $(-A^T\hat z, A\hat x)$ because $\nabla \phi_\mathrm{pd}$ is 
continuous on $\intr(\dom \phi_\mathrm{pd})$.  
Since the operator on right-hand side of~\eqref{e-bcv-fejer-prf} is 
maximal monotone the limit point $(\hat x,\hat z)$ satisfies 
the optimality condition
\[
\begin{bmatrix} -A^T\hat z \\ A\hat x \end{bmatrix} \in 
\begin{bmatrix} \partial f(\hat x)+\nabla h(\hat x) \\
\partial g^\ast(\hat z) \end{bmatrix}
\]
(see \cite[page~27]{Bre:73}, \cite[Lemma~3.2]{Tse:00}).

To show convergence of the entire sequence $(x^{(k)},z^{(k)})$,
we substitute $(\hat x,\hat z)$ in~\eqref{e-bcv-conv-i}:
\[
\cL(x^{(k+1)},\hat z)-\cL(\hat x,z^{(k+1)}) \leq
d(\hat x,\hat z;x^{(k)},z^{(k)})-d(\hat x,\hat z;x^{(k+1)},z^{(k+1)}).
\] 
Since the left-hand side is nonnegative, we have
$d(\hat x,\hat z;x^{(k)},z^{(k)}) \leq 
d(\hat x,\hat z;x^{(k-1)},z^{(k-1)})$
for all $k \geq 1$. This further implies that
\[
d(\hat x,\hat z;x^{(k)},z^{(k)}) \leq 
d(\hat x,\hat z;x^{(k_i)},z^{(k_i)})
\]
for all $k \geq k_i$. By the second additional assumption mentioned above,
the right-hand side converges to zero.
Then the left-hand side also converges to zero and, 
from~(\ref{e-bcv-str-cvx-2})
$x^{(k)} \to \hat x$ and $z^{(k)} \to \hat z$.
\end{proof}

\subsection{Relation to other Bregman proximal algorithms}
\label{s-bcv-connection}

\begin{figure}
\centering
\begin{tikzpicture}[>=latex, font=\small,
                    every text node part/.style={align=center}]
\node (DY) [rectangle,draw]
  {reduced Bregman primal \\ Condat--V\~u};
\node (DR) [rectangle,draw,right=1.5cm of DY]
  {Bregman (primal) \\ Douglas--Rachford};
\node (PG) [rectangle,draw, left=1.5cm of DY]
  {Bregman proximal \\ gradient with shift~\eqref{e-bspg}};
\node (PD3O) [rectangle,draw,above=2.5cm of DY,fill=gray!30]
  {\red{Bregman primal} \\ \red{Condat--V\~u}~\eqref{e-bcv}};
\node (PDHG) [rectangle,draw,above=2.5cm of DR]
  {Bregman \\ (primal) PDHG};
\node (LV) [rectangle,draw,above=2.5cm of PG]
  {Bregman Loris--Verhoeven \\ with shift~\eqref{e-bjv}};
\node (PG0) [rectangle,draw,above=1.0cm of PD3O] 
  {Bregman proximal gradient};
\path (DY.south) -- (DY.south east) coordinate[pos=0.5] (DYe);
\path (DY.south) -- (DY.south west) coordinate[pos=0.5] (DYw);
\path (DR.south) -- (DR.south east) coordinate[pos=0.29] (DRe);
\path (DR.south) -- (DR.south west) coordinate[pos=0.29] (DRw);
\path (PG.south) -- (PG.south east) coordinate[pos=0.5] (PGe);
\path (PG.south) -- (PG.south west) coordinate[pos=0.5] (PGw);
\path (PD3O.north) -- (PD3O.north east) coordinate[pos=0.5] (PD3Oe);
\path (PD3O.north) -- (PD3O.north west) coordinate[pos=0.5] (PD3Ow);
\path (PDHG.north) -- (PDHG.north east) coordinate[pos=0.3] (PDHGe);
\path (PDHG.north) -- (PDHG.north west) coordinate[pos=0.3] (PDHGw);
\path (LV.north) -- (LV.north east) coordinate[pos=0.36] (LVe);
\path (LV.north) -- (LV.north west) coordinate[pos=0.36] (LVw);
\draw[-latex] (DY)   -- node[midway, above] {$h=0$} (DR);
\draw[-latex] (DY)   -- node[midway, above] {$f=0$} (PG);
\draw[-latex] (PD3O) -- node[midway, above] {$h=0$} (PDHG);
\draw[-latex] (PD3O) -- node[midway, above] {$f=0$} (LV);
\draw[latex-] (DY)  -- node[midway, right]  {$A=I$} (PD3O);
\draw[latex-] (DR)  -- node[midway, right]  {$A=I$} (PDHG);
\draw[latex-] (PG)   -- node[midway, right]  {$A=I$} (LV);
\draw[latex-] (PG0)  -- node[midway, right]  {$g=0$} (PD3O);
\end{tikzpicture}
\caption{Proximal algorithms derived from Bregman primal 
Condat--V\~u algorithm~\eqref{e-bcv}.}
\label{t-sum-bcv}
\end{figure}
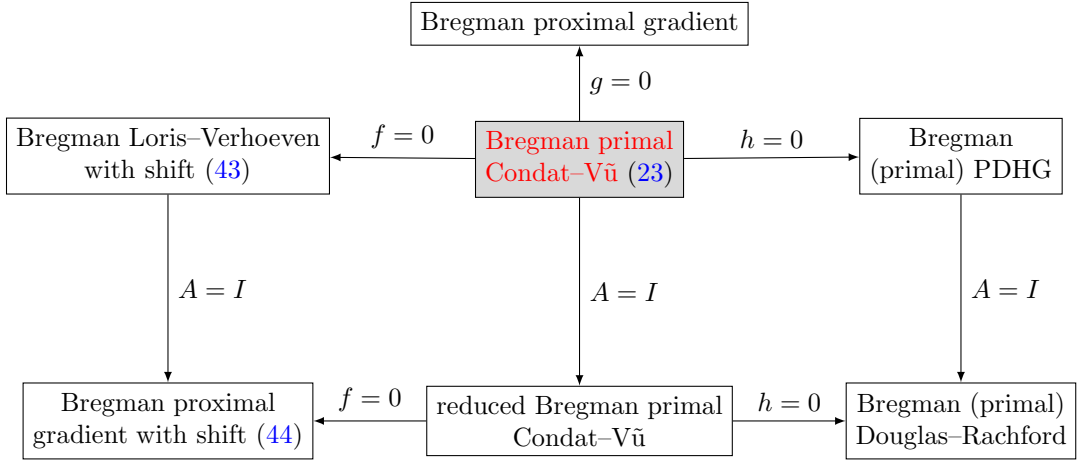

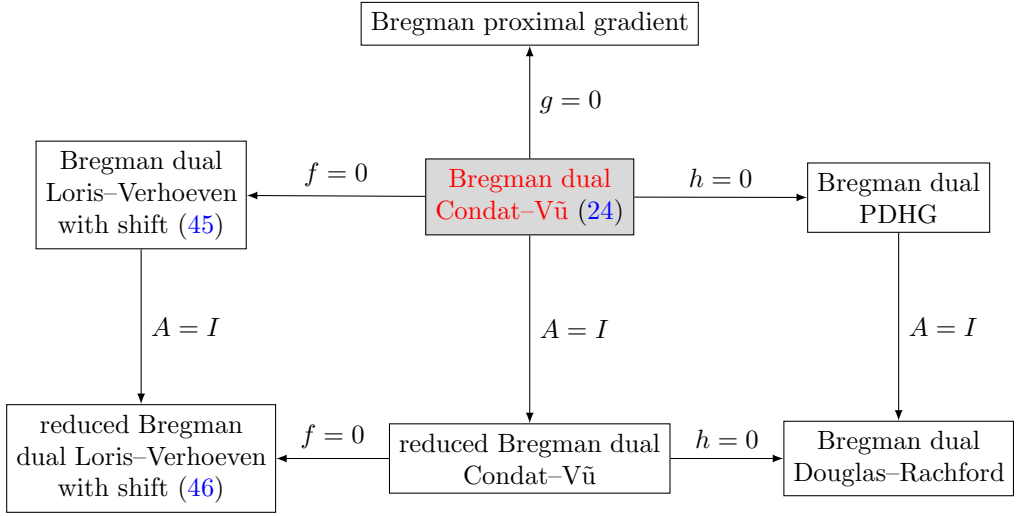
\begin{figure}
\centering
\begin{tikzpicture}[>=latex, font=\small,
                    every text node part/.style={align=center}]
\node (DY) [rectangle,draw]
  {reduced Bregman dual \\ Condat--V\~u};
\node (DR) [rectangle,draw,right=1.5cm of DY]
  {Bregman dual \\ Douglas--Rachford};
\node (PG) [rectangle,draw, left=1.5cm of DY]
  {reduced Bregman \\ dual Loris--Verhoeven \\ 
  with shift~\eqref{e-bspg-dual}};
\node (PD3O) [rectangle,draw,above=2.5cm of DY,fill=gray!30]
  {\red{Bregman dual} \\ \red{Condat--V\~u}~\eqref{e-bcv-dual}};
\node (PDHG) [rectangle,draw,above=2.5cm of DR]
  {Bregman dual \\ PDHG};
\node (LV) [rectangle,draw,above=2.06cm of PG]
  {Bregman dual \\ Loris--Verhoeven \\ with shift~\eqref{e-bjv-dual}};
\node (PG0) [rectangle,draw,above=1.5cm of PD3O] 
  {Bregman proximal gradient};
\path (DY.south) -- (DY.south east) coordinate[pos=0.5] (DYe);
\path (DY.south) -- (DY.south west) coordinate[pos=0.5] (DYw);
\path (DR.south) -- (DR.south east) coordinate[pos=0.29] (DRe);
\path (DR.south) -- (DR.south west) coordinate[pos=0.29] (DRw);
\path (PG.south) -- (PG.south east) coordinate[pos=0.5] (PGe);
\path (PG.south) -- (PG.south west) coordinate[pos=0.5] (PGw);
\path (PD3O.north) -- (PD3O.north east) coordinate[pos=0.5] (PD3Oe);
\path (PD3O.north) -- (PD3O.north west) coordinate[pos=0.5] (PD3Ow);
\path (PDHG.north) -- (PDHG.north east) coordinate[pos=0.3] (PDHGe);
\path (PDHG.north) -- (PDHG.north west) coordinate[pos=0.3] (PDHGw);
\path (LV.north) -- (LV.north east) coordinate[pos=0.36] (LVe);
\path (LV.north) -- (LV.north west) coordinate[pos=0.36] (LVw);
\draw[-latex] (DY)   -- node[midway, above] {$h=0$} (DR);
\draw[-latex] (DY)   -- node[midway, above] {$f=0$} (PG);
\draw[-latex] (PD3O) -- node[midway, above] {$h=0$} (PDHG);
\draw[-latex] (PD3O) -- node[midway, above] {$f=0$} (LV);
\draw[latex-] (DY)  -- node[midway, right]  {$A=I$} (PD3O);
\draw[latex-] (DR)  -- node[midway, right]  {$A=I$} (PDHG);
\draw[latex-] (PG)   -- node[midway, right]  {$A=I$} (LV);
\draw[latex-] (PG0)  -- node[midway, right]  {$g=0$} (PD3O);
\end{tikzpicture}
\caption{Proximal algorithms derived from Bregman dual 
Condat--V\~u algorithm~\eqref{e-bcv-dual}.}
\label{t-sum-bcv-dual}
\end{figure}

Following similar steps as in Section~\ref{s-overview},
we obtain several Bregman proximal splitting methods as special cases
of~\eqref{e-bcv} and~\eqref{e-bcv-dual}.
The connections are summarized in Figure~\ref{t-sum-bcv}
and Figure~\ref{t-sum-bcv-dual}.
A comparison of Figures~\ref{t-sum-cv} and~\ref{t-sum-bcv}
shows that all the reduction relations ($A=I$) are still valid.
However, it is unclear how to apply the ``completion'' operation
to algorithms based on non-Euclidean Bregman distances. 

When $h=0$,~\eqref{e-bcv} reduces to Bregman PDHG~\cite{ChP:16}.
When $g=0$, $g^\ast = \delta_{\{0\}}$ (and assuming $z^{(0)}=0$),
we obtain the Bregman proximal gradient algorithm~\cite{BBT17}.
When $f=0$ in~\eqref{e-bcv}, we obtain the 
\emph{Bregman Loris--Verhoeven algorithm with shift:}
\begin{subequations} \label{e-bjv}
\begin{align}
  x^{(k+1)} &= \argmin_x {\big(
    \inprod{\nabla h(x^{(k)})-A^Tz^{(k)}}{x}+\frac{1}{\tau} 
    d_\primal(x,x^{(k)})\big)} \\
  z^{(k+1)} &= \prox^{\phi_\dual}_{\sigma g^\ast} \big(z^{(k)},
    -\sigma A(2x^{(k+1)}-x^{(k)}) \big).
\end{align}
\end{subequations}
Furthermore, when $A=I$ in~\eqref{e-bcv},
we recover the reduced Bregman primal Condat--V\~u algorithm.
Similarly, setting $A=I$ in Bregman PDHG yields
the Bregman Douglas--Rachford algorithm.
Last, when we set $A=I$ in~\eqref{e-bjv},
we have the \emph{Bregman reduced Loris--Verhoeven algorithm with shift:}
\begin{subequations} \label{e-bspg}
\begin{align}
  x^{(k+1)} &= \argmin_x {\big(\inprod{\nabla h(x^{(k)})-z^{(k)}}{x}
    +\frac{1}{\tau} d_\primal(x,x^{(k)})\big)} \\
  z^{(k+1)} &= \prox^{\phi_\dual}_{\sigma g^\ast}
    \big(z^{(k)},-\sigma (2x^{(k+1)}-x^{(k)}) \big).
\end{align}
\end{subequations}

Similarly, the Bregman dual Condat--V\~u algorithm~\eqref{e-bcv-dual}
can be reduced to some other Bregman proximal splitting methods,
as summarized in Figure~\ref{t-sum-bcv-dual}.
In particular, when $f=0$ in~\eqref{e-bcv-dual}, we obtain
the \emph{Bregman dual Loris--Verhoeven algorithm with shift}:
\begin{subequations} \label{e-bjv-dual}
\begin{align}
z^{(k+1)} &= \prox_{\sigma g^\ast}^{\phi_\dual}
  (z^{(k)}, -\sigma Ax^{(k)}) \\
x^{(k+1)} &= \argmin_x {\big(
  \inprod{A^T(2z^{(k+1)}-z^{(k)})+\nabla h(x^{(k)})}{x}
  + \frac{1}{\tau} d_\primal(x, x^{(k)})\big)}.
\end{align}
\end{subequations}
Moreover, setting $A=I$ in~\eqref{e-bjv-dual} yields the
\emph{reduced Bregman Loris--Verhoeven algorithm with shift:}
\begin{subequations} \label{e-bspg-dual}
\begin{align}
z^{(k+1)} &= \prox_{\sigma g^\ast}^{\phi_\dual}
  (z^{(k)}, -\sigma x^{(k)} ) \\
  x^{(k+1)} &= \argmin_x{\big(
  \inprod{2z^{(k+1)}-z^{(k)} + \nabla h(x^{(k)})}{x}
  + \frac{1}{\tau} d_\primal(x, x^{(k)})\big)}.
\end{align}
\end{subequations}

\section{Bregman dual Condat--V\~u algorithm with line search}
\label{s-ls}

The algorithms~\eqref{e-bcv} and~\eqref{e-bcv-dual} 
use constant parameters $\sigma$ and $\tau$.
The stepsize condition~\eqref{e-bcv-param} involves
the matrix norm $\|A\|$ and the Lipschitz constant~$L$ 
in~\eqref{e-bcv-lip}.
Estimating or bounding $\|A\|$ for a large matrix can be difficult.
As an added complication, the norms $\|\cdot\|_\primal$ and
$\|\cdot\|_\dual$ in the definition of the matrix norm~(\ref{e-norm})
are assumed to be scaled
so that the strong convexity parameters of the primal and dual kernels
are equal to one.  Close bounds on the strong convexity parameters 
may also be difficult to obtain.
Using conservative bounds for $\|A\|$ and $L$ 
results in unnecessarily small values of $\sigma$ and $\tau$, 
and can dramatically slow down the convergence.
Even when the estimates of $\|A\|$ and $L$ are accurate,
the requirements for the stepsizes~\eqref{e-bcv-param} are still
too strict in most iterations, as observed in~\cite{ADH+:21}.
In view of the above arguments,
line search techniques for primal--dual proximal methods have recently 
become an active area of research.
Malitsky and Pock~\cite{MaP:18} proposed a line search technique
for PDHG and the Condat--V\~u algorithm in the Euclidean case.
The algorithm with adaptive parameters in~\cite{VMC21}
focuses on a special case of~\eqref{e-prob} (\ie, $f=0$)
and extends the Loris--Verhoeven algorithm~\eqref{e-lv}.
A Bregman proximal splitting method with line search
is discussed in~\cite{JV22}
and considers the problem~\eqref{e-prob} with $h=0$ and $g=\delta_{\{b\}}$.
In this section, we extend the Bregman dual Condat--V\~u 
algorithm~\eqref{e-bcv-dual} with a varying parameter option,
in which the stepsizes are chosen adaptively without requiring
any estimates or bounds for $\|A\|$ or the strong convexity parameter 
of the kernels.
The algorithm is restricted to problems in the equality constrained
form
\begin{equation} \label{e-lcp-prob}
 \begin{array}{ll}
   \mbox{minimize} & f(x) + h(x) \\
   \mbox{subject to} & Ax=b.
 \end{array}
\end{equation}
This is a special case of~\eqref{e-prob} with $g=\delta_{\{b\}}$,
the indicator function of the singleton $\{b\}$.

The details of the algorithm are discussed in Section~\ref{s-ls-algo}
and a convergence analysis is presented in Section~\ref{s-ls-conv}.
The main conclusion is an $O(1/k)$ rate of ergodic convergence,
consistent with previous results for related algorithms
\cite{MaP:18,JV22}.

\paragraph{Assumptions}
We make the same assumptions as in Section~\ref{s-bcv-ppm}, but define
\[
\phi_\dual (z) = \frac{1}{2} \|z\|^2, \qquad
d_\dual (z,z') = \frac{1}{2} \|z-z'\|^2, \qquad
\|z\|_\dual = \|z\|,
\]
where $\|\cdot\|$ is the Euclidean norm.
The matrix norm $\|A\|$ is defined accordingly as
\[
\|A\|=\sup_{u \neq 0, v \neq 0}
  \frac{\inprod{v}{Au}}{\|v\| \|u\|_\primal}
= \sup_{u \neq 0} \frac{\|Au\|}{\|u\|_\primal}
= \sup_{v \neq 0} \frac{\|A^Tv\|_{\primal,*}}{\|v\|}.
\]

\subsection{Algorithm} \label{s-ls-algo}

The algorithm uses the following iteration,
with starting points $x^{(0)} \in \intr(\dom \phi_\primal) \cap \dom h$
and $z^{(-1)} = z^{(0)}$:
\begin{subequations} \label{e-ls}
\begin{align}
\bar z^{(k+1)} &= z^{(k)} + \theta_k (z^{(k)} - z^{(k-1)}) 
  \label{e-ls-a} \\
x^{(k+1)} &= \prox_{\tau_k f}^{\phi_\primal} \big(
  x^{(k)}, \tau_k (A^T\bar z^{(k+1)}+\nabla h(x^{(k)}))
  \big) \label{e-ls-b} \\
z^{(k+1)} &= z^{(k)} + \sigma_k (Ax^{(k+1)}-b). \label{e-ls-c}
\end{align}
\end{subequations}
With constant parameters $\theta_k=1$, $\sigma_k=\sigma$, $\tau_k = \tau$,
the algorithm can be simplified as
\begin{align*}
x^{(k+1)} &= \prox_{\tau f}^{\phi_\primal} 
  \big( x^{(k)}, \tau A^T(2z^{(k)}-z^{(k-1)})+\tau \nabla h(x^{(k)})) 
  \big) \\
z^{(k+1)} &= z^{(k)} + \sigma (Ax^{(k+1)}-b). 
\end{align*}
Except for the numbering of the dual iterates, this is the
Bregman dual Condat--V\~u algorithm~\eqref{e-bcv-dual} 
applied to~\eqref{e-lcp-prob}.

In the line search algorithm, the parameters 
$\theta_k$, $\tau_k$, $\sigma_k$ are determined by a
backtracking search.
At the start of the algorithm, we set $\tau_{-1}$ and $\sigma_{-1}$ to
some positive values.
To start the search in iteration $k$ we choose
$\bar\theta_k \geq 1$.  For $i=0,1,2,\ldots$,
we set $\theta_k = 2^{-i}\bar\theta_k$,
$\tau_k=\theta_k\tau_{k-1}$, $\sigma_k=\theta_k\sigma_{k-1}$,
and compute $\bar z_{k+1}$, $x_{k+1}$, $z_{k+1}$ using~\eqref{e-ls}.
For some $\delta \in (0,1]$, if
\BEA
& & \inprod{z^{(k+1)}-\bar z^{(k+1)}}{A(x^{(k+1)}-x^{(k)})}
  +h(x^{(k+1)})-h(x^{(k)})-\inprod{\nabla h(x^{(k)})}{x^{(k+1)}-x^{(k)}}
  \nonumber \\
&\leq & \frac{\delta^2}{\tau_k} d_\primal(x^{(k+1)}, x^{(k)})
    + \frac{1}{2\sigma_k} \|\bar z^{(k+1)} - z^{(k+1)}\|^2,
    \label{e-ls-cond}
\EEA
we accept the computed iterates $\bar z^{(k+1)}$, $x^{(k+1)}$, $z^{(k+1)}$
and parameters $\theta_k$, $\sigma_k$, $\tau_k$, and terminate the 
backtracking search.
If~\eqref{e-ls-cond} does not hold, we increment $i$ and continue
the backtracking search.

The backtracking condition~\eqref{e-ls-cond} is similar to
the condition in the line search algorithm for PDHG with Euclidean proximal
operators~\cite[Algorithm~4]{MaP:18}, but it is not identical, even 
in the Euclidean case.
The proposed condition is weaker and allows larger stepsizes
than the condition in~\cite[Algorithm~4]{MaP:18}.

\subsection{Convergence analysis} \label{s-ls-conv}

The proof strategy is the same as in~\cite[Section~3.3]{JV22},
extended to account for the function $h$.
The main conclusion is an $O(1/k)$ rate of ergodic convergence,
shown in equation~\eqref{e-ls-conv-1/k}.

\subsubsection{Lower bound on algorithm parameters}

We first show that the stepsizes are bounded below by
\BEQ \label{e-ls-param-lbnd}
  \tau_k \geq \tau_\mathrm{min} \triangleq \min{\Big\{\tau_{-1},
  \frac{-L+\sqrt{L^2+4\delta^2\beta\|A\|^2}}{4\beta\|A\|^2}\Big\}},
  \qquad
  \sigma_k \geq \sigma_\mathrm{min} \triangleq \beta\tau_\mathrm{min},
\EEQ
where $\beta=\sigma_{-1}/\tau_{-1}$.
The lower bounds imply that the backtracking eventually terminates
with positive stepsizes $\sigma_k$ and $\tau_k$.
\begin{proof}
Applying~\eqref{e-bcv-delta} with $\tau=\tau_k$, $\sigma=\sigma_k$, $\delta_1=\delta^2-\tau_k L$ and $\delta_2=1$,
together with the Lipschitz condition~\eqref{e-bcv-lip},
we see that the backtracking condition~\eqref{e-ls-cond} holds
at iteration $k$ if $0 < \delta < 1$ and
\[
\tau_k \sigma_k\|A\|^2 + \tau_k L \leq \delta^2.
\]
Then mathematical induction can be used to prove~\eqref{e-ls-param-lbnd}.
The two lower bounds~\eqref{e-ls-param-lbnd} hold at $k=0$
by the definition of $\tau_\mathrm{min}$ and $\sigma_\mathrm{min}$.
Now assume $\tau_{k-1} \geq \tau_\mathrm{min}$,
$\sigma_{k-1} \geq \sigma_\mathrm{min}$,
and consider the $k$th iteration.
The first attempt of $\theta_k$ is $\theta_k=\bar \theta_k \geq 1$.
If this value is accepted, then
\[
\tau_k = \bar\theta_k \tau_{k-1} \geq \tau_{k-1} \geq \tau_\mathrm{min},
\qquad
\sigma_k = \bar\theta_k \sigma_{k-1} \geq \sigma_{k-1} \geq
  \sigma_\mathrm{min}.
\]
Otherwise, one or more backtracking steps are needed.
Denote by $\tilde \theta_k$ the last rejected value.
Then $\tilde\theta_k^2\tau_{k-1}^2\beta\|A\|^2 + 
\tilde\theta_k \tau_{k-1} L >\delta^2$
and the accepted $\theta_k$ satisfies
\[
\theta_k = \frac{\tilde \theta_k}{2} \geq
\frac{-L+\sqrt{L^2+4\delta^2\beta\|A\|^2}}%
  {4\tau_{k-1}\beta\|A\|^2}.
\]
Therefore,
\[
\tau_k=\theta_k\tau_{k-1} >
\frac{-L+\sqrt{L^2+4\delta^2\beta\|A\|^2}}{4\beta\|A\|^2}, \qquad
\sigma_k=\beta\tau_k \geq \beta\tau_\mathrm{min}. \qedhere
\]
\end{proof}

\subsubsection{One-iteration analysis} \label{e-cv-1iter}

The iterates $x^{(k+1)}$, $z^{(k+1)}$, $\bar z^{(k+1)}$
generated by the algorithm~\eqref{e-ls} satisfy
\begin{align}
\cL(x^{(k+1)},z)-\cL(x, \bar z^{(k+1)})
  &\leq \frac{1}{\tau_k} 
  \left(d_\primal(x,x^{(k)}) - d_\primal(x,x^{(k+1)}) -
  (1-\delta^2) d_\primal(x^{(k+1)},x^{(k)})\right) \nonumber \\
&\phantom{=} + \frac{1}{2\sigma_k} \left(
  \|z-z^{(k)}\|^2 - \|z-z^{(k+1)}\|^2 - 
  \|\bar z^{(k+1)} - z^{(k)}\|^2 \right) \label{e-ls-conv-i}
\end{align}
for all $x \in \dom f \cap \dom \phi_\primal$ and all $z$.
Here $\cL(x,z)=f(x)+h(x)+\inprod{z}{Ax-b}$.  
\begin{proof}
The optimality condition for the primal prox-operator~\eqref{e-ls-b}
gives
\BEAS
\lefteqn{f(x^{(k+1)})-f(x)} \\
&\leq & \frac{1}{\tau_k}
  \big(d_\primal(x,x^{(k)})-d_\primal(x,x^{(k+1)})-
  d_\primal(x^{(k+1)}, x^{(k)})\big)
  +\inprod{A^T\bar z^{(k+1)}+\nabla h(x^{(k)})}{x-x^{(k+1)}}
\EEAS
for all $x\in \dom f \cap \dom\phi_\primal$.
Hence
\BEA
\lefteqn{f(x^{(k+1)}) + h(x^{(k+1)}) - f(x) - h(x)} \nonumber \\
  &\leq & \frac{1}{\tau_k} (d_\primal(x,x^{(k)}) - 
  d_\primal(x,x^{(k+1)}) - d_\primal(x^{(k+1)}, x^{(k)}))
  + \inprod{A^T\bar z_{k+1}}{x - x^{(k+1)}} \nonumber \\
& &\mbox{} + h(x^{(k+1)})-h(x)+\inprod{\nabla h(x^{(k)})}{x-x^{(k+1)}}
  \nonumber \\
&\leq & \frac{1}{\tau_k} (d_\primal(x,x^{(k)}) - 
 d_\primal(x,x^{(k+1)}) - d_\primal(x^{(k+1)}, x^{(k)}))
  + \inprod{A^T\bar z^{(k+1)}}{x - x^{(k+1)}} \nonumber \\
& &\mbox{} +h(x^{(k+1)})-h(x^{(k)})
  -\inprod{\nabla h(x^{(k)})}{x^{(k+1)}-x^{(k)}}.
  \label{e-ls-conv-1}
\EEA
The second inequality follows from the convexity of $h$, \ie,
$h(x) \geq h(x^{(k)}) + \inprod{\nabla h(x^{(k)})}{x-x^{(k)}}$.
The dual update~\eqref{e-ls-c} implies that
\BEQ \label{e-ls-conv-2}
\inprod{z-z^{(k+1)}}{Ax^{(k+1)}-b}
= \frac{1}{\sigma_k} \inprod{z-z^{(k+1)}}{z^{(k+1)}-z^{(k)}}
\quad \mbox{for all $z$.}
\EEQ
This equality at $k=i-1$ is
\begin{align}
\inprod{z-z^{(i)}}{Ax^{(i)}-b}
  &= \frac{1}{\sigma_{i-1}} \inprod{z-z^{(i)}}{z^{(i)}-z^{(i-1)}} 
  \nonumber \\
&= \frac{1}{2\sigma_{i-1}} \left( \|z-z^{(i-1)}\|^2
  - \|z-z^{(i)}\|^2 - \|z^{(i)} - z^{(i-1)} \|^2 \right).
  \label{e-ls-conv-3}
\end{align}
The equality~\eqref{e-ls-conv-2} at $k=i-2$ is
\begin{align*}
\inprod{z-z^{(i-1)}}{Ax^{(i-1)}-b}
  &= \frac{1}{\sigma_{i-2}} \inprod{z-z^{(i-1)}}{z^{(i-1)}-z^{(i-2)}} \\
&= \frac{\theta_{i-1}}{\sigma_{i-1}}
  \inprod{z-z^{(i-1)}}{z^{(i-1)}-z^{(i-2)}} \\
&= \frac{1}{\sigma_{i-1}} \inprod{z-z^{(i-1)}}{\bar z^{(i)}-z^{(i-1)}}.
\end{align*}
We evaluate this at $z=z^{(i)}$ and add it to the equality
at $z=z^{(i-2)}$ multiplied by $\theta_{i-1}$:
\begin{align}
\inprod{z^{(i)}-\bar z^{(i)}}{Ax^{(i-1)}-b}
  &= \frac{1}{\sigma_{i-1}}
  \inprod{z^{(i)}-\bar z^{(i)}}{\bar z^{(i)} - z^{(i-1)}} \nonumber \\
&= \frac{1}{2\sigma_{i-1}} \left(\|z^{(i)}-z^{(i-1)}\|^2
  -\|z^{(i)}-\bar z^{(i)}\|^2-\|\bar z^{(i)}-z^{(i-1)}\|^2 \right).
  \label{e-ls-conv-4}
\end{align}
Now we combine~\eqref{e-ls-conv-1} for $k=i-1$,
with~\eqref{e-ls-conv-3} and~\eqref{e-ls-conv-4}. For $i \geq 1$,
\BEAS
\lefteqn{\cL(x^{(i)},z) - \cL(x, \bar z^{(i)})} \nonumber \\
&=& f(x^{(i)})+h(x^{(i)})+\inprod{z}{Ax^{(i)}-b}
  -f(x)-h(x)-\inprod{\bar z^{(i)}}{Ax-b} \nonumber \\
&\leq & \frac{1}{\tau_{i-1}}
  \Big(d_\primal(x,x^{(i-1)})-d_\primal(x,x^{(i)})
  - d_\primal(x^{(i)},x^{(i-1)})\Big)
  + \inprod{A^T\bar z^{(i)}}{x-x^{(i)}} + \inprod{z}{Ax^{(i)}-b} \nonumber\\
& &\mbox{} - \inprod{\bar z^{(i)}}{Ax-b} + h(x^{(i)}) - h(x^{(i-1)})
  -\inprod{\nabla h(x^{(i-1)})}{x^{(i)}-x^{(i-1)}} \nonumber \\
&=& \frac{1}{\tau_{i-1}}
  \Big(d_\primal(x,x^{(i-1)}) - d_\primal(x,x^{(i)}) 
  - d_\primal(x^{(i)},x^{(i-1)})\Big)
  + \inprod{z-\bar z^{(i)}}{Ax^{(i)}-b} \nonumber \\
& &\mbox{} +h(x^{(i)})-h(x^{(i-1)})
  -\inprod{\nabla h(x^{(i-1)})}{x^{(i)}-x^{(i-1)}} \nonumber \\
&=& \frac{1}{\tau_{i-1}}
  \Big(d_\primal(x,x^{(i-1)}) - d_\primal(x,x^{(i)})
  - d_\primal(x^{(i)},x^{(i-1)})\Big) \nonumber \\
& & \mbox{} +\inprod{z^{(i)}-\bar z^{(i)}}{A(x^{(i)}-x^{(i-1)})}
  + \inprod{z-z^{(i)}}{Ax^{(i)}-b}
  + \inprod{z^{(i)}-\bar z^{(i)}}{Ax^{(i-1)}-b} \nonumber \\
& &\mbox{} +h(x^{(i)})-h(x^{(i-1)})
  -\inprod{\nabla h(x^{(i-1)})}{x^{(i)}-x^{(i-1)}} \nonumber \\
&=& \frac{1}{\tau_{i-1}}
  \Big(d_\primal(x,x^{(i-1)}) - d_\primal(x,x^{(i)})
  - d_\primal(x^{(i)},x^{(i-1)}))\Big) \nonumber \\
& &\mbox{} +\frac{1}{2\sigma_{i-1}}
  \Big(\|z-z^{(i-1)}\|^2-\|z-z^{(i)}\|^2-\|\bar z^{(i)}-z^{(i-1)}\|^2
  -\|\bar z^{(i)}-z^{(i)}\|^2\Big) \nonumber \\
& &\mbox{}+\inprod{A^T(z^{(i)}-\bar z^{(i)})}{x^{(i)}-x^{(i-1)}}
  +h(x^{(i)})-h(x^{(i-1)})
  -\inprod{\nabla h(x^{(i-1)})}{x^{(i)}-x^{(i-1)})} \\
&\leq & \frac{1}{\tau_{i-1}} 
  \left(d_\primal(x,x^{(i-1)}) - d_\primal(x,x^{(i)}) -
  (1-\delta^2) d_\primal(x^{(i)},x^{(i-1)})\right) \nonumber \\
& &\mbox{} + \frac{1}{2\sigma_{i-1}} \left(
  \|z-z^{(i-1)}\|^2 - \|z-z^{(i)}\|^2 - 
  \|\bar z^{(i)} - z^{(i-1)}\|^2 \right),
\EEAS
which is the desired result~\eqref{e-ls-conv-i}.
The first inequality follows from~\eqref{e-ls-conv-1}.
In the second last step we substitute~\eqref{e-ls-conv-3}
and~\eqref{e-ls-conv-4}.
The last step uses the line search exit condition~\eqref{e-ls-cond}
at $k=i-1$.
\end{proof}

\subsubsection{Ergodic convergence}

We define the averaged primal and dual sequences
\[
x^{(k)}_\avg = \frac{1}{\sum_{i=1}^k \tau_{i-1}}
\sum_{i=1}^k \tau_{i-1}x^{(i)}, \qquad
\bar z^{(k)}_\avg = \frac{1}{\sum_{i=1}^k \tau_{i-1}}
\sum_{i=1}^k \tau_{i-1} \bar z^{(i)}
\]
for $k \geq 1$.
We show that
\begin{align}
\cL(x^{(k)}_\avg,z) - \cL(x,\bar z^{(k)}_\avg)
  &\leq \frac{1}{\sum_{i=1}^k \tau_{i-1}}
  \big(d_\primal(x,x^{(0)}) + \frac{1}{2\beta} \|z-z^{(0)}\|^2\big)
  \label{e-ls-conv} \\ 
&\leq \frac{1}{k \tau_\mathrm{min}}
  \big(d_\primal(x,x^{(0)}) + \frac{1}{2\beta} \|z-z^{(0)}\|^2\big)
  \label{e-ls-conv-1/k}
\end{align}
for all $x \in \dom f \cap \dom \phi_\primal$ and all $z$.
This holds for any choice of $\delta \in (0,1]$ in~\eqref{e-ls-cond}.
If we compare~\eqref{e-ls-conv} and~\eqref{e-bcv-conv}, we note that the 
two left-hand sides involve different dual iterates 
($\bar z_\avg^{(k)}$ as opposed to $z_\avg^{(k)}$).

\begin{proof}
From~\eqref{e-ls-conv-i},
\[
\cL(x^{(i)},z)-\cL(x, \bar z^{(i)}) \leq \frac{1}{\tau_{i-1}}
\Big(d_\primal(x,x^{(i-1)}) - d_\primal(x,x^{(i)})
+ \frac{1}{2\beta} \|z-z^{(i-1)}\|^2
- \frac{1}{2\beta} \|z-z^{(i)}\|^2 \Big).
\]
Since $\cL(u,v)$ is convex in $u$ and affine in $v$,
\BEA
\lefteqn{(\sum_{i=1}^k \tau_{i-1}) 
  \big(\cL(x^{(k)}_\avg, z)-\cL(x, \bar z^{(k)}_\avg)\big)} \nonumber \\
&\leq &\sum_{i=1}^k \tau_{i-1} (\cL(x^{(i)}, z)-\cL(x, \bar z^{(i)})) 
  \nonumber \\
&\leq & d_\primal(x,x^{(0)}) - d_\primal(x,x^{(k)}) +
  \frac{1}{2\beta} (\|z-z^{(0)}\|^2 - \|z-z^{(k)}\|^2) \nonumber \\
&\leq & d_\primal(x,x^{(0)}) + \frac{1}{2\beta} \|z-z^{(0)}\|^2.
  \label{e-ls-conv-k}
\EEA
Dividing by $\sum_{i=1}^k \tau_{i-1}$ gives~\eqref{e-ls-conv}.
\end{proof}

Substituting $x=x^\star$ and $z=z^\star$ in~\eqref{e-ls-conv-k} yields
\[
f(x^{(k)}_\avg) + h(x^{(k)}_\avg) + \inprod{z}{Ax^{(k)}_\avg-b} 
- f(x^\star) - h(z^\star) \leq
\frac{1}{\sum_{i=1}^k\tau_{i-1}} \big(d_\primal(x^\star,x^{(0)})
+ \frac{1}{2\beta} \|z^\star-z^{(0)}\|^2 \big),
\]
since $Ax^\star=b$.
More generally, suppose $X \subseteq \dom f \cap \dom \phi_\primal$
is a compact convex set containing an optimal solution $x^\star$
in its interior,
and $Z = \{z \mid \|z\| \leq \gamma\}$ contains a dual optimal $z^\star$,
then the merit function~$\eta$ defined in~\eqref{e-merit} satisfies
\begin{align*}
\eta(x_\avg^{(k)}, \bar z_\avg^{(k)})
&= \sup_{z \in Z} \mathcal L(x_\avg^{(k)},z)
  - \inf_{x \in X} \mathcal L(x,\bar z_\avg^{(k)}) \\
&\leq \frac{1}{\sum_{i=1}^k\tau_{i-1}} \Big(
  \sup_{x\in X} d_\primal(x,x^{(0)})
  + \frac{1}{2\beta} (\gamma + \|z^{(0)}\|)^2 \Big) \\
&\leq \frac{1}{k \tau_\mathrm{min}}
  \Big(\sup_{x \in X} d_\primal(x,x^{(0)})
  + \frac{1}{2\beta} (\gamma + \|z^{(0)}\|)^2 \Big).
\end{align*}
The second line follows from~\eqref{e-ls-conv}
and the third line follows from~\eqref{e-ls-param-lbnd}.

\subsubsection{Monotonicity properties and convergence of iterates}

For $x=x^\star$, $z=z^\star$, the left-hand side of~\eqref{e-ls-conv-i}
is nonnegative and we obtain
\BEAS
\lefteqn{d_\primal(x^\star,x^{(k+1)})
  +\frac{1}{2\beta}\|z^\star-z^{(k+1)}\|^2} \\
&\leq & d_\primal(x^\star,x^{(k)})+\frac{1}{2\beta}\|z^\star-z^{(k)}\|^2
  -\big((1-\delta^2)d_\primal(x^{(k+1)},x^{(k)})
  +\frac{1}{2\beta}\|z^{(k+1)}-z^{(k)}\|^2 \big) \\
&\leq & d_\primal (x^\star,x^{(k)})+\frac{1}{2\beta}\|z^\star-z^{(k)}\|^2
\EEAS
for $k \geq 0$. Moreover,
\[
\sum_{i=0}^k \Big( (1-\delta^2) (d_\primal(x^{(i+1)},x^{(i)})
  +\frac{1}{2\beta} \|\bar z^{(i+1)}-z^{(i)} \|^2 \Big) 
\leq d_\primal (x^\star, x^{(0)}) 
  +\frac{1}{2\beta} \|z^\star-\bar z^{(0)}\|^2.
\]
These inequalities hold for any value $\delta \in (0,1]$.
In particular, the last inequality implies that
$\bar z^{(i+1)}-z^{(i)} \to 0$.
When $\delta < 1$ it also implies that
$d_\primal(x^{(i+1)}, x^{(i)}) \to 0$ and,
by the strong convexity assumption on $\phi_\primal$,
that $x^{(i+1)}-x^{(i)} \to 0$.
With additional assumptions similar to those in Section~\ref{s-fejer},
one can show the convergence of iterates;
see~\cite[Section~3.3.4]{JV22}.

\section{Bregman PD3O algorithm} \label{s-bpd3o}

\begin{figure}
\centering
\begin{tikzpicture}
\node at (4,-0.3) {\small {$1/L$}};
\fill[fill=gray!20!white] plot[domain=0.28:4, smooth] (\x,{1/\x-1/4)}) -- 
    plot[domain=4:0, smooth] (\x,{0}) --
    plot[domain=0:0.2, smooth] (\x, {95/28});
\fill[fill=gray!50!white] plot[domain=0.28:4, smooth] (\x,{1/\x-1/4)}) -- 
    plot[domain=4:0.32, smooth] (\x,{1/\x+1/4});
\filldraw[blue] (4,0) circle(1.2pt);
\draw[-latex] (0, 0) -- (5, 0) node[right] {$\tau$};
\draw[-latex] (0, 0) -- (0, 3.5) node[above] {$\sigma$};
\draw[domain=0.28:4, smooth, thick, variable=\x, blue]
  plot ({\x}, {1/\x-1/4});
\draw[domain=0.32:4, smooth, thick, variable=\x, red]  
  plot ({\x}, {1/\x+1/4});
\end{tikzpicture}
\caption{Acceptable stepsizes in Condat--V\~u algorithms and PD3O.
We assume the same matrix norm $\|A\|$ and Lipschitz constant $L$
are used in the analysis of the two algorithms.
The light gray region under the blue curve is defined by the 
inequality for the Condat--V\~u algorithms in~\eqref{e-param}.
The region under the red curve shows the values allowed 
by the stepsized conditions for PD3O.
}
\label{f-param-1}
\end{figure}
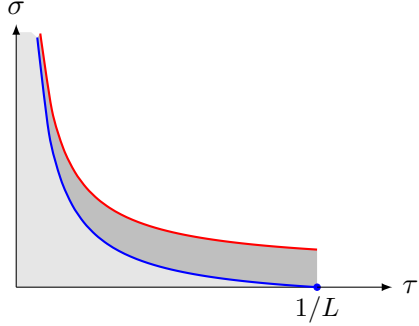

In this section we propose the \emph{Bregman PD3O algorithm},
another Bregman proximal method for the problem~\eqref{e-prob}.
Bregman PD3O also involves two generalized distances,
$d_\primal$ and $d_\dual$, generated by $\phi_\primal$ and $\phi_\dual$,
respectively, and it consists of the iterations
\begin{subequations} \label{e-bpd3o}
\begin{align}
  x^{(k+1)} &= \prox^{\phi_\primal}_{\tau f}
    (x^{(k)},\tau A^Tz^{(k)}+\tau\nabla h(x^{(k)})) \\
  z^{(k+1)} &= \prox^{\phi_\dual}_{\sigma g^\ast} (z^{(k)},
    -\sigma A(2x^{(k+1)}-x^{(k)}
    +\tau (\nabla h(x^{(k)})-\nabla h(x^{(k+1)})) ) ).
\end{align}
\end{subequations}
The only difference between Bregman PD3O and Bregman primal Condat--V\~u 
algorithm~\eqref{e-bcv} is the additional term
$\tau (\nabla h(x^{(k)})-\nabla h(x^{(k+1)}) )$.
Thus the two algorithms~\eqref{e-bcv} and~\eqref{e-bpd3o} reduce to
the same method when $h$ is absent from problem~\eqref{e-prob}.
The additional term allows PD3O to use larger stepsizes than 
the Condat--V\~u algorithm.
If we use the same matrix norm $\|A\|$ and Lipschitz constant $L$ in
the analysis for the two methods, then the conditions are
\begin{equation} \label{e-param}
  \begin{array}{ll}
    \mbox{Condat--V\~u:} & \sigma \tau \|A\|^2 + \tau L \leq 1 \\
    \mbox{PD3O:} & \sigma \tau \|A\|^2 \leq 1, \;\; \tau \leq 1/L.
  \end{array}
\end{equation}
The range of possible parameters is illustrated in Figure~\ref{f-param-1}.

In Section~\ref{s-bpd3o-conv} we provide the detailed convergence analysis
of the Bregman PD3O method.
The connections between Bregman PD3O and several other Bregman proximal
methods are discussed in Section~\ref{s-bpd3o-connection}.

\paragraph{Assumptions}
Throughout Section~\ref{s-bpd3o} we make the following assumptions.
The kernel functions $\phi_\primal$ and $\phi_\dual$ are 
$1$-strongly convex with respect to the Euclidean norm and
an arbitrary norm $\|\cdot\|_\dual$, respectively:
\begin{equation} \label{e-bpd3o-str-cvx}
d_\primal(x,x^\prime) \geq \frac{1}{2}\|x-x^\prime\|^2, \qquad
d_\dual(z,z^\prime) \geq \frac{1}{2}\|z-z^\prime\|_\dual^2.
\end{equation} 
The assumptions that the strong convexity constants are one can
be made without loss of generality, by scaling the distances.
The definition of $\|A\|$ follows~\eqref{e-norm} and reduces to
\[
\|A\|=\sup_{v\neq 0} \frac{\|Av\|_{\dual,\ast}}{\|v\|}.
\] 
We also assume that the gradient of $h$ is $L$-Lipschitz continuous
with respect to the Euclidean norm: $\dom h = \reals^n$ and
\BEQ \label{e-bpd3o-lip-1}
h(y)-h(x)-\inprod{\nabla h(x)}{y-x} \leq \frac{L}{2}\|y-x\|^2,
\quad \mbox{for any } x,y \in \dom h.
\EEQ
The parameters $\tau$ and $\sigma$ must satisfy
\begin{equation} \label{e-bpd3o-param}
\sigma \tau \|A\|^2 \leq 1, \qquad \tau \leq 1/L.
\end{equation} 
Finally, we assume that the optimality condition~\eqref{e-opt-cond} has a
solution $(x^\star,z^\star) \in \dom \phi_\primal \times \dom \phi_\dual$.

Note that~\eqref{e-bpd3o-lip-1} is a stronger assumption
than~\eqref{e-bcv-lip}.
(Combined with the first inequality in~\eqref{e-bpd3o-str-cvx},
it implies~\eqref{e-bcv-lip}.)
We will use the following consequence of~\eqref{e-bpd3o-lip-1}:
\begin{equation} \label{e-bpd3o-lip-2}
h(y)-h(x)-\inprod{\nabla h(x)}{y-x} 
\geq \frac{1}{2L} \|\nabla h(y)-\nabla h(x)\|^2
\end{equation} 
for all $x,y$ \cite[Theorem 2.1.5]{Nesterov18}.

\subsection{Convergence analysis} \label{s-bpd3o-conv}

\subsubsection{A primal--dual Bregman distance}
We introduce a primal--dual kernel
\[
\phi_\mathrm{pd3o}(x,y,z) = \frac{1}{\tau} \phi_\primal(x)
  +\frac{1}{\sigma} \phi_\dual(z)+\frac{\tau}{2} \|y\|^2 
  -\inprod{y}{x}-\inprod{z}{A(x-\tau y)},
\]
where $\sigma,\tau > 0$.  If $\phi_\mathrm{pd3o}$ is convex,
the generated Bregman distance is given by
\begin{align}
d_\mathrm{pd3o}(x,y,z;x^\prime,y^\prime,z^\prime)
  &= \frac{1}{\tau} d_\primal(x,x^\prime)
  +\frac{1}{\sigma} d_\dual(z,z^\prime)
  +\frac{\tau}{2} \|y-y^\prime\|^2 \nonumber \\
&\phantom{=} - \inprod{y-y^\prime}{x-x^\prime}
  -\inprod{z-z^\prime}{A(x-x^\prime)}
  +\tau \inprod{z-z^\prime}{A(y-y^\prime)}. \label{e-bpd3o-dist}
\end{align}
We now show that $\phi_\mathrm{pd3o}$ is convex if
$\sigma \tau \|A\|^2 \leq 1$.
\begin{proof}
It is sufficient to show that $d_\mathrm{pd3o}$ is nonnegative:
\begin{align}
d_\mathrm{pd3o}(x,y,z;x^\prime,y^\prime,z^\prime)
  &\geq \frac{1}{2\tau} \|x-x^\prime\|^2
  +\frac{\tau}{2} \|A^T(z-z^\prime)\|^2
  +\frac{\tau}{2} \|y-y^\prime\|^2 \nonumber \\
&\phantom{=} - \inprod{y-y^\prime}{x-x^\prime}
  -\inprod{z-z^\prime}{A(x-x^\prime)}
  +\tau \inprod{z-z^\prime}{A(y-y^\prime)} \nonumber \\
&= \frac{1}{2} \Big\|\frac{1}{\sqrt\tau}(x-x^\prime)
  -\sqrt\tau (y-y^\prime)
  -\sqrt\tau A^T(z-z^\prime) \Big\|^2 \nonumber \\
&\geq 0. \label{e-bpd3o-norm}
\end{align}
In step 1 we use the strong convexity assumption~\eqref{e-bpd3o-str-cvx},
the definition of $\|A\|$~\eqref{e-norm}
with $\|\cdot\|_\primal=\|\cdot\|$,
and the assumption $\sigma\tau\|A\|^2 \leq 1$. 
The bound on $d_\dual(z,z')$ follows from
\[
\frac{1}{\sigma}d_\dual(z,z^\prime)
\geq \frac{1}{2\sigma}\|z-z^\prime\|_\dual^2
\geq \frac{\|A^T(z-z^\prime)\|^2}{2\sigma\|A\|^2}
\geq \frac{\tau}{2} \|A^T(z-z^\prime)\|^2.
\qedhere
\]
\end{proof}

Note that the convexity of $\phi_\mathrm{pd3o}$ only requires the first
inequality in the stepsize condition~\eqref{e-bpd3o-param}.
Although the Bregman PD3O algorithm~\eqref{e-bpd3o} is not the
Bregman proximal point method for the Bregman kernel $\phi_\mathrm{pd3o}$,
the distance $d_\mathrm{pd3o}$ will appear in the key 
inequality~\eqref{e-bpd3o-conv-i} of the convergence analysis. 

\subsubsection{One-iteration analysis}

We first show that the iterates $x^{(k+1)}$, $z^{(k+1)}$ generated by
Bregman PD3O~\eqref{e-bpd3o} satisfy
\BEA
\lefteqn{\cL(x^{(k+1)},z)-\cL(x,z^{(k+1)})} \nonumber \\
&\leq & d_\mathrm{pd3o}
  \big(x,\nabla h(x), z;x^{(k)},\nabla h(x^{(k)}),z^{(k)}\big)
  -d_\mathrm{pd3o} \big(
    x,\nabla h(x), z;x^{(k+1)},\nabla h(x^{(k+1)}),z^{(k+1)}
  \big) \nonumber \\
& &\mbox{} -d_\mathrm{pd3o}\big(x^{(k+1)},\nabla h(x), z^{(k+1)};%
    x^{(k)},\nabla h(x^{(k)}),z^{(k)} \big) \label{e-bpd3o-conv-i}
\EEA
for all $x \in \dom f \cap \dom \phi_\primal$ and
$z \in \dom g^\ast \cap \dom \phi_\dual$.
\begin{proof}
Recall that Bregman PD3O differs from the Bregman primal Condat--V\~u
algorithm~\eqref{e-bcv} only in an additional term
in the dual update.
The proof in Section~\ref{e-cv-1iter} therefore applies up 
to~(\ref{e-bcv-conv-prf-1}), with
\[
\tilde x = 2x^{(k+1)} - x^{(k)} + \tau (\nabla h(x^{(k)}) - 
 \nabla h(x^{(k+1)})), \qquad \tilde z = z^{(k)}.
\]
Substituting the above $(\tilde x, \tilde z)$ into~\eqref{e-bcv-conv-prf-2}
and applying the definition of~$d_-$~\eqref{e-bcv-dminus} yields
\BEAS
\lefteqn{\cL(x^{(k+1)},z)-\cL(x,z^{(k+1)})} \\
&\leq & d_-(x,z;x^{(k)},z^{(k)}) - d_-(x,z;x^{(k+1)},z^{(k+1)})
  - d_-(x^{(k+1)},z^{(k+1)};x^{(k)},z^{(k)}) \nonumber \\
& &\mbox{} -
 \tau \inprod{A^T(z-z^{(k+1)})}{\nabla h(x^{(k)})-\nabla h(x^{(k+1)})} \\
& &\mbox{} +h(x^{(k+1)})-h(x)
  +\inprod{\nabla h(x^{(k)})}{x-x^{(k+1)}} \nonumber \\
&=& d_-(x,z;x^{(k)},z^{(k)})+\frac{\tau}{2}
  \|\nabla h(x)-\nabla h(x^{(k)})\|^2 \\
& &\mbox{} -\inprod{(x-\tau A^Tz)%
  -(x^{(k)}-\tau A^Tz^{(k)})}%
  {\nabla h(x)-\nabla h(x^{(k)})} \\
& &\mbox{} -\Big(d_-(x,z;x^{(k+1)},z^{(k+1)})+\frac{\tau}{2}
  \|\nabla h(x)-\nabla h(x^{(k+1)})\|^2 \\
& &\mbox{} -\big\langle{x-\tau A^Tz
  -(x^{(k+1)}-\tau A^Tz^{(k+1)} )},%
  {\nabla h(x)-\nabla h(x^{(k+1)})} \big\rangle \Big) \\
& &\mbox{} -\Big(d_-(x^{(k+1)}, z^{(k+1)}; x^{(k)}, z^{(k)})
  +\frac{\tau}{2}
  \|\nabla h(x)-\nabla h(x^{(k)})\|^2 \\
& &\mbox{} -\big\langle
  {(x^{(k+1)}-\tau A^Tz^{(k+1)}) -(x^{(k)}-\tau A^Tz^{(k)})},%
  {\nabla h(x)-\nabla h(x^{(k)})} \big\rangle \Big) \\
& &\mbox{} -(
  h(x)-h(x^{(k+1)})-\inprod{\nabla h(x^{(k+1)})}{x-x^{(k+1)}}
  -\frac{\tau}{2}\|\nabla h(x)-\nabla h(x^{(k+1)})\|^2) \\
&=& d_\mathrm{pd3o}(x,\nabla h(x),z;x^{(k)},\nabla h(x^{(k)}),z^{(k)})
  -d_\mathrm{pd3o}
    (x,\nabla h(x),z;x^{(k+1)},\nabla h(x^{(k+1)}),z^{(k+1)}) \\
& &\mbox{} -d_\mathrm{pd3o} (x^{(k+1)},\nabla h(x),z^{(k+1)};%
    x^{(k)},\nabla h(x^{(k)}),z^{(k)}) \\
& &\mbox{} -(
  h(x)-h(x^{(k+1)})-\inprod{\nabla h(x^{(k+1)})}{x-x^{(k+1)}}
  -\frac{\tau}{2}\|\nabla h(x)-\nabla h(x^{(k+1)})\|^2) \\
&\leq & d_\mathrm{pd3o}
    (x,\nabla h(x),z;x^{(k)},\nabla h(x^{(k)}),z^{(k)})
  -d_\mathrm{pd3o}
    (x,\nabla h(x),z;x^{(k+1)},\nabla h(x^{(k+1)}),z^{(k+1)}) \\
& &\mbox{} -d_\mathrm{pd3o} (x^{(k+1)},\nabla h(x),z^{(k+1)};%
    x^{(k)},\nabla h(x^{(k)}),z^{(k)}) \\
&\leq & d_\mathrm{pd3o}
    (x,\nabla h(x),z;x^{(k)},\nabla h(x^{(k)}),z^{(k)})
  -d_\mathrm{pd3o}
    (x,\nabla h(x),z;x^{(k+1)},\nabla h(x^{(k+1)}),z^{(k+1)}).
\EEAS
Step 3 follows from definition of~$d_\mathrm{pd3o}$~\eqref{e-bpd3o-dist}.
In step 4 we use the Lipschitz condition~\eqref{e-bpd3o-lip-2}
and the second inequality in the stepsize condition~\eqref{e-bpd3o-param}.
The last step follows from the fact that 
$d_\mathrm{pd3o}$ is nonnegative~\eqref{e-bpd3o-norm}.
\end{proof}

\subsubsection{Ergodic convergence}

The iterates generated by Bregman PD3O~\eqref{e-bpd3o} satisfy
\[
\cL(x^{(k)}_\avg, z)-\cL(x, z^{(k)}_\avg)
\leq \frac{3}{k} \Big(\frac{2}{\tau} d_\primal (x,x^{(0)})
+ \frac{1}{\sigma} d_\dual (z,z^{(0)}) \Big),
\]
for all $x \in \dom f \cap \dom \phi_\primal$ and 
all $z \in \dom g^\ast \cap \dom \phi_\dual$,
where the averaged iterates are defined in~\eqref{e-bcv-avg}.
\begin{proof}
From~\eqref{e-bpd3o-conv-i},
since $\cL(u,v)$ is convex in $u$ and concave in $v$,
\BEAS
\lefteqn{\cL(x^\avg_k, z) - \cL(x, z^\avg_k)} \\
  &\leq &\frac{1}{k} \sum_{i=1}^k
  \big(\cL(x_i,z)-\cL(x,z_i)\big) \\
&\leq & \frac{1}{k}
  d_\mathrm{pd3o}(x,\nabla h(x), z;x_0,\nabla h(x_0), z_0) \\
&\leq & \frac{3}{k} \Big(\frac{1}{\tau} d_\primal (x,x_0)
  + \frac{1}{\sigma} d_\dual (z,z_0) + \frac{\tau}{2}
  \|\nabla h(x)-\nabla h(x_0)\|^2 \Big) \\
&\leq & \frac{3}{k} \Big(\frac{1}{\tau} d_\primal (x,x_0)
  + \frac{1}{\sigma} d_\dual (z,z_0) + \frac{\tau L^2}{2}
  \|x-x_0\|^2 \Big) \\
&\leq & \frac{3}{k} \Big(\frac{2}{\tau} d_\primal (x,x_0)
  + \frac{1}{\sigma} d_\dual (z,z_0) \Big)
\EEAS
for all $x\in \dom f \cap \dom \phi_\primal$ and
$z \in \dom g^\ast \cap \dom \phi_\dual$.
The third inequality follows from~(\ref{e-bpd3o-dist}):
\begin{align*}
&d_\mathrm{pd3o}(x,y,z;x',y',z') \\
&\mbox{} \leq \frac{1}{\tau}d_\mathrm p(x,x')
 + \frac{1}{\sigma} d_\mathrm d(z,z') + \frac{\tau}{2} \|y-y'\|^2
 + \|y-y'\|\|x-x'\| \\
& \mbox{} \phantom{\leq}
 + \|A\| \|x-x'\|\|z-z'\|_\mathrm d + \|A\| \|y-y'\|\|z-z'\|_\mathrm d \\
&\mbox{} \leq \frac{1}{\tau}d_\mathrm p(x,x')
 + \frac{1}{\sigma} d_\mathrm d(z,z') + \frac{\tau}{2} \|y-y'\|^2
 + \frac{1}{2\tau} \|x-x'\|^2 + \frac{\tau}{2} \|y-y'\|^2 \\
& \mbox{} \phantom{\leq}
 + \frac{1}{2\tau} \|x-x'\|^2 + \frac{1}{2\sigma} \|z-z'\|_\mathrm d^2
 + \frac{1}{2\tau} \|y-y'\|^2 + \frac{1}{2\sigma} \|z-z'\|_\mathrm d^2 \\
&\mbox{} \leq \frac{3}{\tau} d_\primal(x,x')
 + \frac{3}{\sigma} d_\dual(z,z') + \frac{3\tau}{2} \|y-y'\|^2.
\qedhere
\end{align*}
\end{proof}

\subsection{Relation to other Bregman proximal algorithms}
\label{s-bpd3o-connection}

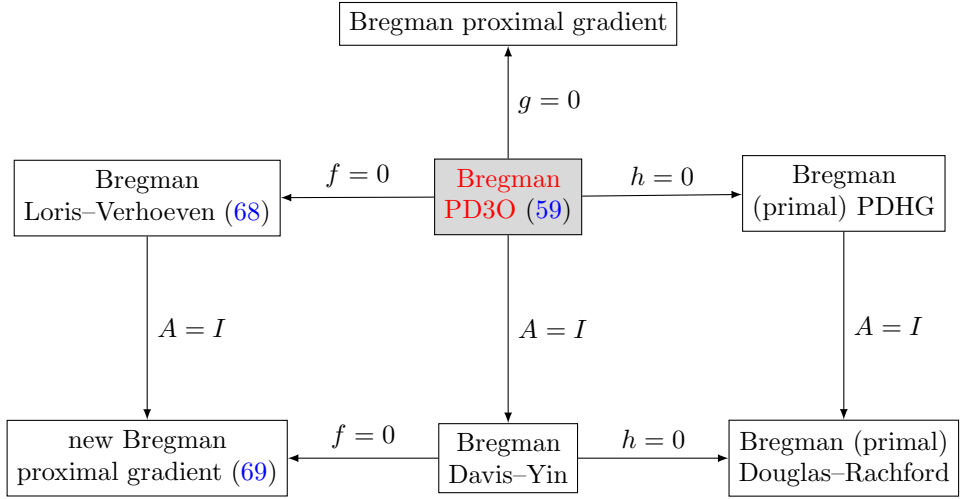
\begin{figure}
\centering
\begin{tikzpicture}[>=latex, font=\small,
                    every text node part/.style={align=center}]
\node (DY) [rectangle,draw] {Bregman \\ Davis--Yin};
\node (DR) [rectangle,draw,right=2cm of DY]
  {Bregman (primal) \\ Douglas--Rachford};
\node (PG) [rectangle,draw, left=2cm of DY]
  {new Bregman \\ proximal gradient~\eqref{e-bpg}};
\node (PD3O) [rectangle,draw,above=2.5cm of DY,fill=gray!30]
  {\red{Bregman} \\ \red{PD3O}~\eqref{e-bpd3o}};
\node (PDHG) [rectangle,draw,above=2.5cm of DR]
  {Bregman \\ (primal) PDHG};
\node (LV) [rectangle,draw,above=2.44cm of PG]
  {Bregman \\ Loris--Verhoeven~\eqref{e-blv}};
\node (PG0) [rectangle,draw,above=1.5cm of PD3O]
  {Bregman proximal gradient};
\path (DY.south) -- (DY.south east) coordinate[pos=0.33] (DYe);
\path (DY.south) -- (DY.south west) coordinate[pos=0.33] (DYw);
\path (DR.south) -- (DR.south east) coordinate[pos=0.25] (DRe);
\path (DR.south) -- (DR.south west) coordinate[pos=0.24] (DRw);
\path (PG.south) -- (PG.south east) coordinate[pos=0.5] (PGe);
\path (PG.south) -- (PG.south west) coordinate[pos=0.5] (PGw);
\path (PD3O.north) -- (PD3O.north east) coordinate[pos=0.44] (PD3Oe);
\path (PD3O.north) -- (PD3O.north west) coordinate[pos=0.44] (PD3Ow);
\path (PDHG.north) -- (PDHG.north east) coordinate[pos=0.26] (PDHGe);
\path (PDHG.north) -- (PDHG.north west) coordinate[pos=0.26] (PDHGw);
\path (LV.north) -- (LV.north east) coordinate[pos=0.36] (LVe);
\path (LV.north) -- (LV.north west) coordinate[pos=0.36] (LVw);
\draw[-latex] (DY)   -- node[midway, above] {$h=0$} (DR);
\draw[-latex] (DY)   -- node[midway, above] {$f=0$} (PG);
\draw[-latex] (PD3O) -- node[midway, above] {$h=0$} (PDHG);
\draw[-latex] (PD3O) -- node[midway, above] {$f=0$} (LV);
\draw[latex-] (DY)  -- node[midway, right]  {$A=I$} (PD3O);
\draw[latex-] (DR)  -- node[midway, right]  {$A=I$} (PDHG);
\draw[latex-] (PG)  -- node[midway, right]  {$A=I$} (LV);
\draw[latex-] (PG0)  -- node[midway, right]  {$g=0$} (PD3O);
\end{tikzpicture}
\caption{Proximal algorithms derived from Bregman PD3O.}
\label{t-sum-bpd3o}
\end{figure}

The proposed algorithm~\eqref{e-bpd3o} can be viewed as
an extension to PD3O~\eqref{e-pd3o} using generalized distances,
and reduces to several Bregman proximal methods by reduction.
These algorithms can also be organized into a diagram similar to
Figure~\ref{t-sum-pd3o}.
Figure~\ref{t-sum-bpd3o} starts from Bregman PD3O~\eqref{e-bpd3o},
and summarizes its connection to several Bregman proximal methods.
When $h=0$,~\eqref{e-bpd3o} reduces to Bregman PDHG,
and when $g=0$,~\eqref{e-bpd3o} reduces to the Bregman proximal
gradient algorithm.
The Bregman Loris--Verhoeven algorithm is Bregman PD3O with $f=0$:
\begin{subequations} \label{e-blv}
\begin{align}
  x^{(k+1)} &= \argmin_x {\big(
    \inprod{\nabla h(x^{(k)})-A^Tz^{(k)}}{x}+\frac{1}{\tau}
    d_\primal (x,x^{(k)}) \big)} \\
  z^{(k+1)} &= \prox^{\phi_\dual}_{\sigma g^\ast} \big(z^{(k)},
    -\sigma A\big(2x^{(k+1)}-x^{(k)} + \tau
    \big(\nabla h(x^{(k)})-\nabla h(x^{(k+1)})\big)\big) \big).
\end{align}
\end{subequations}
This algorithm has been discussed in~\cite{CST21} under the name NEPAPC.
Setting $A=I$ (with $\sigma=1/\tau$),
we obtain a new variant of Bregman proximal gradient algorithm:
\begin{subequations} \label{e-bpg}
\begin{align}
  x^{(k+1)} &= \argmin_x {\big(
    \inprod{\nabla h(x^{(k)})-z^{(k)}}{x}+\frac{1}{\tau}
    d_\primal (x,x^{(k)}) \big)} \\
  z^{(k+1)} &= \prox^{\phi_\dual}_{\tau^{-1} g^\ast} \big(z^{(k)},
    -\frac{1}{\tau}A\big(2x^{(k+1)}-x^{(k)})
    -A(\nabla h(x^{(k)})-\nabla h(x^{(k+1)})) \big).
\end{align}
\end{subequations}
The difference between~\eqref{e-bpg} and~\eqref{e-bspg} is
the additional term $\tau(\nabla h(x^{(k)})-\nabla h(x^{(k+1)}))$,
the same as the difference between~\eqref{e-bcv} and~\eqref{e-bpd3o}.
When the Euclidean proximal operator is used,
\eqref{e-bpg} reduces to the proximal gradient method.
However, the new algorithm~\eqref{e-bpg} does not seem
to be equivalent to the Bregman proximal gradient algorithm
due to the lack of Moreau decomposition in the generalized case.
Nevertheless, the new algorithm~\eqref{e-bpg}
may still be interesting on its own,
especially when the generalized proximal operator of $g^\ast$ is easy to
compute while the (Euclidean or generalized) proximal operator of $g$
is computationally expensive.
Finally, setting $A=I$ (and $\sigma=1/\tau$)
in Bregman PD3O~\eqref{e-bpd3o} gives a Bregman Davis--Yin algorithm.

\section{Numerical experiment} \label{s-exp}

In this section we evaluate the performance of the Bregman primal
Condat--V\~u algorithm~\eqref{e-bcv}, 
Bregman dual Condat--V\~u algorithm with line search~\eqref{e-ls},
and Bregman PD3O~\eqref{e-bpd3o}.
The main goal of the example is to validate and illustrate
the difference in the stepsize conditions~\eqref{e-param}, and
the usefulness of the line search procedure.
We consider the convex optimization problem
\begin{equation} \label{e-exp-prob}
\begin{array}{ll}
  \mbox{minimize} & \psi(x) = \lambda \|Ax\|_1+\tfrac{1}{2}\|Cx-b\|^2 \\
  \mbox{subject to} & \ones^Tx=1, \quad x \succeq 0,
\end{array}
\end{equation} 
where $x \in \reals^n$ is the optimization variable,
$C \in \reals^{m \times n}$,
and $A \in \reals^{(n-1) \times n}$ is the difference matrix
\BEQ \label{e-diff-mat}
A = \begin{bmatrix}[r]
  -1 & 1 & 0 & \cdots & 0 & 0 \\
  0 & -1 & 1 & \cdots & 0 & 0 \\
  \vdots & \vdots & \vdots & & \vdots & \vdots \\
  0 & 0 & 0 & \cdots & -1 & 1
\end{bmatrix}.
\EEQ 
This problem is of the form of~\eqref{e-prob} with
\[
f(x)=\delta_H(x), \qquad g(y)=\lambda \|y\|_1, \qquad
g^\ast(z) = \begin{cases}
  0 \quad &\|z\|_\infty \leq \lambda \\
  +\infty &\mbox{otherwise,}
\end{cases} \qquad h(x)=\frac{1}{2}\|Cx-b\|^2,
\] 
and $\delta_H$ is the indicator function of the hyperplane
$H=\{x \in \reals^n \mid \ones^Tx=1\}$.
We use the relative entropy distance 
\[
d_\primal(x,y) = \sum_{i=1}^n(x_i \log(x_i/y_i) - x_i + y_i),
\qquad \dom d_\primal = \reals^n_+ \times \reals^n_{++}.
\]
in the primal space.
This distance is 1-strongly convex with respect to
$\ell_1$-norm~\cite{BeT:09b} (and also $\ell_2$-norm).
With the relative entropy distance, all the primal iterates $x^{(k)}$
remain feasible.
In the dual space we use the Euclidean distance.
Thus, the matrix norm~\eqref{e-norm} in the stepsize 
condition~\eqref{e-bcv-param} for the Bregman Condat--V\~u
algorithms is the (1,2)-operator norm
\[
\|A\|_{1,2}=\sup_{v \neq 0} \frac{\|Av\|}{\|v\|_1}
=\max_{i=1,\ldots,n} \|a_i\|=\sqrt 2,
\] 
where $a_i$ is the $i$th column of $A$.
In the Bregman PD3O algorithm, we use the squared Euclidean
distance $d_\primal(x,y) = \tfrac{1}{2} \|x-y\|^2$,
and the matrix norm in the stepsize condition~(\ref{e-bpd3o-param})
is the spectral norm $\|A\|_2$.
For the difference matrix~(\ref{e-diff-mat}),
$\|A\|_2$ is bounded above by $2$, and very close to this upper bound
for large $n$.  

The Lipschitz constant for $h$ with respect to the $\ell_1$-norm is
the largest absolute value of the elements in $C^TC$,
\ie, $L_1=\max_{i,j} |(C^TC)_{ij}|$.   This value is used in 
the stepsize condition~\eqref{e-bcv-param} for the
Bregman Condat--V\~u algorithms.
The Lipschitz constant with respect to the $\ell_2$-norm is
$L_2=\|C\|_2^2$, which is used in the stepsize
condition~(\ref{e-bpd3o-param}) for Bregman PD3O.

The matrix norms and Lipschitz constants are summarized as follows:
\[
\begin{array}{lcc}
& \mbox{matrix norm} & \mbox{Lipschitz constant} \\
\mbox{Bregman Condat--V\~u} \quad & \|A\|_{1,2}=\sqrt 2 \quad &
  L_1=\max_{i,j}|(C^TC)_{ij}| \\
\mbox{Bregman PD3O} & \|A\|_2 \leq 2 & L_2=\|C\|_2^2.
\end{array}
\] 
In the example we use the exact values of $L_1$ and $L_2$, 

The Bregman proximal operator of $f$ has a closed-form solution:
\[
\prox_f^\phi(y,a) = \frac{1}{\sum_{i=1}^n y_i e^{-a_i}} \begin{bmatrix}
 y_1 e^{-a_1} \\ \vdots \\ y_n e^{-a_n}
\end{bmatrix},
\]
and the (Euclidean) proximal operator of $g^\ast$ is 
the projection onto the infinity norm ball:
\[
\prox_{g^\ast}(z)_i = \begin{cases}
  \lambda \quad & z_i>\lambda \\
  z_i & |z_i| \leq \lambda \\
  -\lambda & z_i<-\lambda.
\end{cases}
\] 

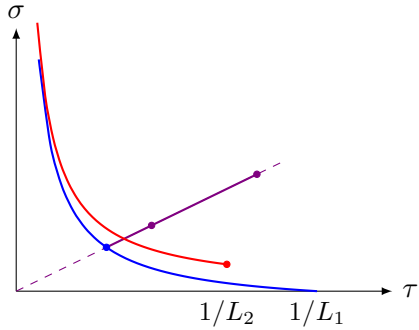
\begin{figure}
\centering
\begin{tikzpicture}
\node at (4,-0.3) {\small {$1/L_1$}};
\node at (14/5,-0.3) {\small {$1/L_2$}};
\draw[-latex] (0, 0) -- (5, 0) node[right] {$\tau$};
\draw[-latex] (0, 0) -- (0, 3.5) node[above] {$\sigma$};
\draw[dashed, violet] (0,0) -- (18/5, 7/4);
\draw[thick, violet] (6/5,7/12) -- (16/5, 14/9);
\draw[domain=0.30:4, smooth, thick, variable=\x, blue]
  plot ({\x}, {1/\x-1/4});
\draw[domain=0.28:2.8, smooth, thick, variable=\x, red]  
  plot ({\x}, {1/\x});
\filldraw[red] (14/5,5/14) circle(1.2pt);
\filldraw[blue] (6/5,7/12) circle(1.2pt);
\filldraw[violet] (9/5,7/8) circle(1.2pt);
\filldraw[violet] (16/5,14/9) circle(1.2pt);
\end{tikzpicture}
\caption{The blue and red curves show the boundaries of the stepsize 
regions
for Bregman Condat--V\~u algorithms and Bregman PD3O, respectively.
The blue and red points indicate the chosen parameters in~\eqref{e-setup}
(red for for PD3O, blue for Condat--V\~u).
In the Bregman dual Condat--V\~u algorithm with line search,
the stepsizes are selected on the dashed straight line.  
The solid line segment shows the range of 
stepsizes that were selected, with dots indicating the largest, median, 
and smallest stepsizes.}
\label{f-param-2}
\end{figure}

The experiment is carried out in Python 3.6 on a desktop
with an Intel Core i5 2.4GHz CPU and 8GB RAM.
We set $m=500$ and $n=10,000$.
The elements in the matrix $C \in \reals^{m \times n}$
and $b \in \reals^m$ are randomly generated 
from independent standard Gaussian distributions.
For the constant stepsize option, we choose
\begin{equation} \label{e-setup}
\begin{array}{lll}
\mbox{Condat-V\~u} & \sigma=L_1/2 & \tau=1/(2L_1) \\
\mbox{PD3O} & \sigma=L_2/4 & \tau=1/L_2.
\end{array}
\end{equation} 
These two choices, as well as the range of possible parameters,
are illustrated in Figure~\ref{f-param-2}.
The two choices are on the blue and red curve, respectively,
and satisfy the requirement~\eqref{e-param} with equality.
For the line search algorithm, we set $\bar \theta_k=1.2$ to 
encourage more aggressive updates, and $\beta=\sigma_{-1}/\tau_{-1}=L_1^2$,
which is consistent with the choice in~\eqref{e-setup}.

\begin{figure}
\begin{minipage}{0.48\textwidth}
\centering
\begin{tikzpicture}[font=\scriptsize, scale = 0.95]
\centering
\begin{axis}[xmin=-5, xmax=70, ymode=log,
  xlabel = {time (sec)},
  ylabel = {${(\psi(x^{(k)})-\psi^\star)}/{\psi^\star}$},
  xlabel near ticks, ylabel near ticks, yticklabel style={rotate=90},
  legend pos=south west, legend cell align={left},
  label style={font=\footnotesize}, tick label style={font=\tiny}]
\addplot[color=blue] table [x=cvtime, y=cvobj, col sep=comma]%
  {result-1.csv};
\addplot[color=violet] table [x=lstime, y=lsobj, col sep=comma]%
  {result-1.csv};
\addplot[color=red] table [x=pdtime, y=pdobj, col sep=comma]%
  {result-1.csv};
\legend{Bregman CV, Bregman CV w. LS, Bregman PD3O};
\end{axis}
\end{tikzpicture}
\end{minipage}%
\hfill
\begin{minipage}{0.48\textwidth}
\centering
\begin{tikzpicture}[font=\scriptsize, scale = 0.95]
\centering
\begin{axis}[xmin=-500, xmax=22000, ymode=log,
  xlabel = {number of iterations},
  ylabel = {${(\psi(x^{(k)})-\psi^\star)}/{\psi^\star}$},
  xlabel near ticks, ylabel near ticks, yticklabel style={rotate=90},
  legend pos=south west, legend cell align={left},
  label style={font=\footnotesize}, tick label style={font=\tiny}]
\addplot[color=blue] table [x=cviter, y=cvobj, col sep=comma]%
  {result-1.csv};
\addplot[color=violet] table [x=lsiter, y=lsobj, col sep=comma]%
  {result-1.csv};
\addplot[color=red] table [x=pditer, y=pdobj, col sep=comma]%
  {result-1.csv};
\legend{Bregman CV, Bregman CV w. LS, Bregman PD3O};
\end{axis}
\end{tikzpicture}
\end{minipage}

\vspace{6pt}
\begin{minipage}{0.48\textwidth}
\centering
\begin{tikzpicture}[font=\scriptsize, scale = 0.95]
\centering
\begin{axis}[xmin=-5, xmax=70, ymode=log,
  xlabel = {time (sec)},
  ylabel = {${(\psi(x^{(k)})-\psi^\star)}/{\psi^\star}$},
  xlabel near ticks, ylabel near ticks, yticklabel style={rotate=90},
  legend pos=south west, legend cell align={left},
  label style={font=\footnotesize}, tick label style={font=\tiny}]
\addplot[color=blue] table [x=cvtime, y=cvobj, col sep=comma]%
  {result-2.csv};
\addplot[color=violet] table [x=lstime, y=lsobj, col sep=comma]%
  {result-2.csv};
\addplot[color=red] table [x=pdtime, y=pdobj, col sep=comma]%
  {result-2.csv};
\legend{Bregman CV, Bregman CV w. LS, Bregman PD3O};
\end{axis}
\end{tikzpicture}
\end{minipage}%
\hfill
\begin{minipage}{0.48\textwidth}
\centering
\begin{tikzpicture}[font=\scriptsize, scale = 0.95]
\centering
\begin{axis}[xmin=-500, xmax=8000, ymode=log,
  xlabel = {number of iterations},
  ylabel = {${(\psi(x^{(k)})-\psi^\star)}/{\psi^\star}$},
  xlabel near ticks, ylabel near ticks, yticklabel style={rotate=90},
  legend pos=south west, legend cell align={left},
  label style={font=\footnotesize}, tick label style={font=\tiny}]
\addplot[color=blue] table [x=cviter, y=cvobj, col sep=comma]%
  {result-2.csv};
\addplot[color=violet] table [x=lsiter, y=lsobj, col sep=comma]%
  {result-2.csv};
\addplot[color=red] table [x=pditer, y=pdobj, col sep=comma]%
  {result-2.csv};
\legend{Bregman CV, Bregman CV w. LS, Bregman PD3O};
\end{axis}
\end{tikzpicture}
\end{minipage}
\caption{Comparison of three algorithms (Bregman primal Condat--V\~u,
Bregman dual Condat--V\~u with line search, and Bregman PD3O)
in terms of objective values.
The top two figures plot the relative error of the function value
versus CPU time and number of iterations
for one problem instance~\eqref{e-exp-prob}, respectively.
The bottom two figures correspond to another problem instance.
}
\label{f-result}
\end{figure}
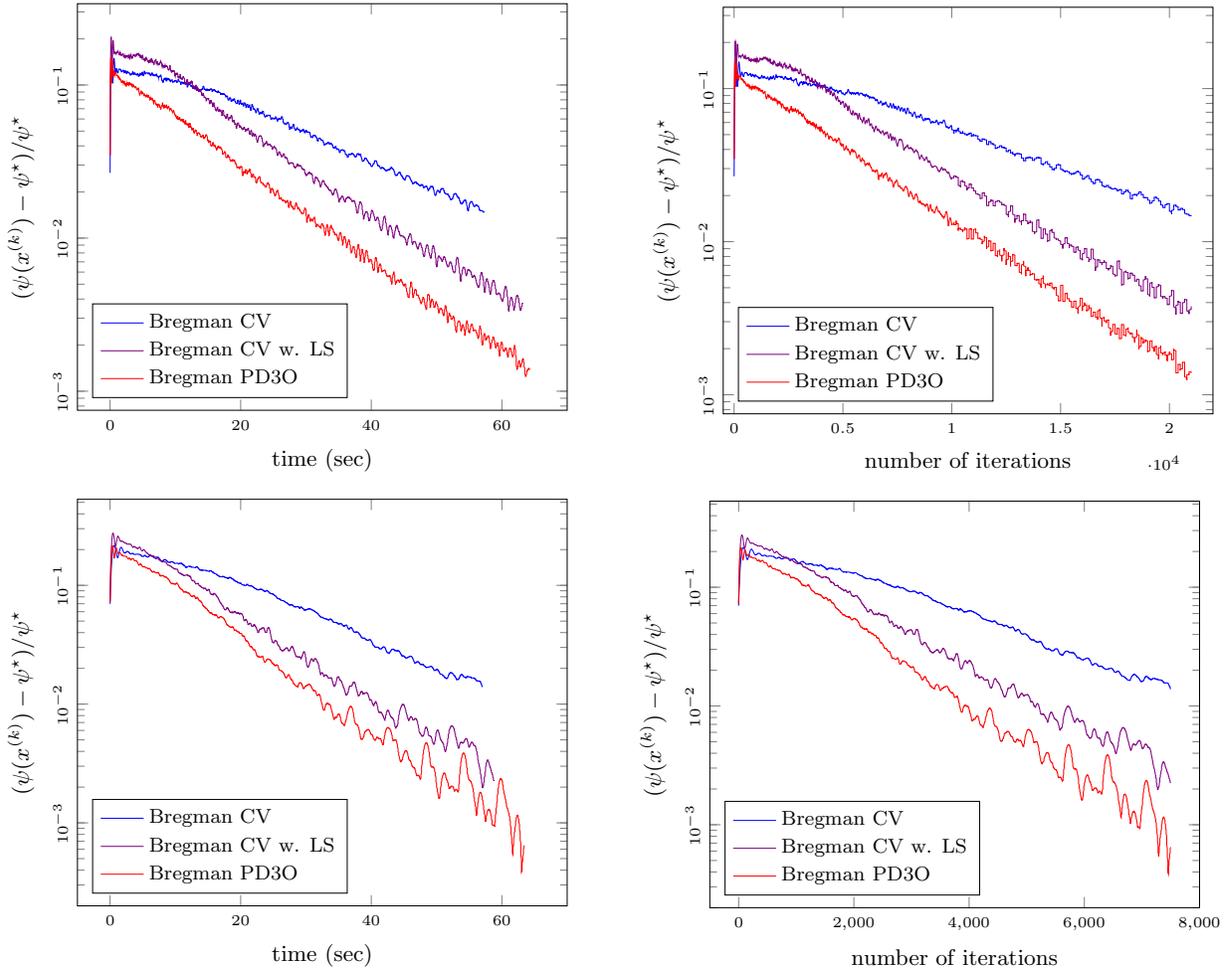

We solve the problem~\eqref{e-exp-prob} using 
the Bregman primal Condat--V\~u algorithm~\eqref{e-bcv},
the Bregman dual Condat--V\~u algorithm with line search~\eqref{e-ls},
and Bregman PD3O~\eqref{e-bpd3o}.
Figure~\ref{f-result} reports the relative distance
between the function values to the optimal value $\psi^\star$,
which is computed via CVXPY~\cite{DCB:14}.
Comparison between the Bregman primal Condat--V\~u algorithm
and Bregman PD3O shows that Bregman PD3O converges faster.
Figure~\ref{f-result} also compares the performance
between the Bregman primal Condat--V\~u algorithm
with constant stepsizes and Bregman dual algorithm with line search.
One can see clearly that the line search significantly improves the 
convergence.
On the other hand, the line search does not add much computation
overhead, as the plots of the CPU time and the number of iterations
are roughly identical.
In these experiments Bregman PD3O
and the Bregman dual Condat--V\~u algorithm with line search have
a similar performance,
without one algorithm being conclusively better than the other.

\section{Conclusions}

We presented two variants of Bregman Condat--V\~u algorithms,
introduced a line search technique
for the Bregman dual Condat--V\~u algorithm for equality-constrained
problems, and proposed a Bregman extension to PD3O.
Many open questions remain.
It is unclear how to use Bregman distances in PDDY,
and how to extend the line search technique to Bregman PD3O,
the Bregman primal Condat--V\~u algorithm, 
and the more general problem~\eqref{e-prob}.
Moreover, in the current backtracking technique the ratio of the
primal and dual stepsizes is fixed.
A further improvement would be to relax this constraint
\cite{MaP:18,ADH+:21}.


\end{document}